\documentclass[notitlepage,11pt]{article}

\usepackage{amssymb,amsfonts,amsmath,geometry,esint,comment,upgreek,cases,graphicx,bbm}
\usepackage{mathtools}

\usepackage{color}

\catcode`\@=11
\@addtoreset{equation}{section}

\catcode`\@=12

\allowdisplaybreaks

\def\Ds{\left(-\Delta\right)^{\!s}\!} 
\def\Dshalf{\left(-\Delta\right)^{\!\frac {s}{2}}\!}  

\newcommand{\bDelta}{\Delta\!\!\!\!\scriptstyle{\Delta}\displaystyle}
\newcommand{\eltaD}{\reflectbox{$\Delta$}}
\newcommand{\bbDelta}{\scalebox{0.7}{$\eltaD$}\hspace{0.2pt}\!\!\!\!\!{\eltaD}\!\!\displaystyle}
\newcommand{\bnabla}{\hspace{2.5pt} \reflectbox{ {\rotatebox[origin=c]{180}{$\bbDelta$}}}}

\def\R{\mathbb{R}}

\def\N{\mathbb{N}}

\def\f{\varphi}
\def\e{\varepsilon}
\def\P{\mathrm{P}}
\def\E{\mathbb E}
\def\D{\mathcal D}

\def\me{\mathrm{\bf e}}

\def\irn{\int\limits_{\R^n}}
\def\iirn{\iint\limits_{\R^{n+1}}}
\def\div{{\rm div}}
\def\nn{|\!|\!|}

\def\proof{\noindent{\textbf{Proof. }}}
\def\QED{\hfill {$\square$}\goodbreak \medskip}

\newtheorem{Theorem}{Theorem}[section]
\newtheorem{Lemma}[Theorem]{Lemma}

\newtheorem{Corollary}[Theorem]{Corollary}
\newtheorem{Remark}[Theorem]{Remark}
\newtheorem{Definition}[Theorem]{Definition}

\begin{document}

\title{The $s$-polyharmonic extension problem\\ and higher-order fractional Laplacians}

\author{Gabriele Cora\footnote{DMIF, Universit\`a di Udine, Italy.
Email: {gabriele.cora@uniud.it}, orcid.org/0000-0002-0090-5470.
}
~ and Roberta Musina\footnote{
DMIF, Universit\`a di Udine, Italy.
Email: {roberta.musina@uniud.it}, orcid.org/0000-0003-4835-8004.}}

\date{}

\maketitle

\noindent
{\footnotesize {\bf Abstract.}
We provide a detailed description of the relationships between the fractional Laplacian of order $2s\in(0,n)$ on $\R^n$ and the
{\em $s$-polyharmonic} extension operator. 
}


\section{Introduction}

The seminal paper  \cite{CS} by Caffarelli and Silvestre marked a turning point in the study of the fractional Laplacian
$\Ds$ for $s\in(0,1)$. It was shown that for any sufficiently regular function $u$ on $\R^n$, there exists a unique solution $\E_s[u]$ to the 
Dirichlet problem
$$
-\div(y^{1-2s}\nabla U)=0\quad\text{in $\R^{n+1}_+$}~,\qquad  U(~\!\cdot~\!,0)=u\,~\!,
$$
which satisfies
$$ 
\iint\limits_{\R^{n+1}_+}y^{1-2s}|\nabla\E_s[u]|^2~\!dz=d_s\irn|\Dshalf u|^2~\!dx\, ,\qquad
-\lim_{y \to 0^+}y^{1-2s}\partial_y \E_s[u] =  d_s\,(-\Delta)^s u~\!.
$$
Here $z=(x,y)\in \R^{n+1}_+\equiv\R^n\times(0,\infty)$ and 
$d_s= \frac{1}{\Gamma(s)}~\!2^{1-2s}\Gamma(1- s)$.
Further, the {\em $s$-harmonic extension} $\E_s[u]$
of $u$ 
can be expressed via convolution with a Poisson kernel, 
$$
\E_{s}[u](x,y) = (u* \P_{\!s}^y)(x)~, \quad\text{where}\quad \P^y_s(x)=\frac{\Gamma\big(\frac{n}{2}+s\big)}{\pi^{\frac{n}{2}}\Gamma(s)} \frac{y^{2s}}{(|x|^2+y^2)^\frac{n+2s}{2}}~,
$$
and can be easily characterized as the unique solution to the minimization problem
$$
\inf_{V \in \mathcal D^{1;1-2s}(\R^{n+1}_+) \atop V(\cdot ,0)=u}
\iint\limits_{\R^{n+1}_+}
y^{1-2s}|\nabla V|^2~\!dz~\!,
$$
where $\mathcal D^{1;1-2s}(\R^{n+1}_+)$ is a suitably defined energy space.

Attempts to fully extend the Caffarelli-Silvestre approach to higher orders have already been made,
via conformal geometry techniques, starting with \cite{CdMG} (see also the related papers \cite{CC, MMG, GR, JX, RS} and references therein\footnote
{\footnotesize{the unpublished preprint \cite{Ya} (partially included also in \cite{CY}) contains several inaccuracies.}}). For  more recent and comprehensive results, we refer to \cite{Case}. We cite also  \cite{FF}, where
$s\in(1,2)$ is assumed and a different approach is used. 

It has already been  observed that for any $s$ having integer part $[s]\in[1,s)$, the fractional Laplacian $\Ds$  can be recovered by solving a
(possibly) degenerate
elliptic boundary value problem of order $2(1+[s])$, and that the extension operator $\E_s$ can still be used to represent its solution. On the other hand, as far as we know, a clean and complete description of the natural functional framework for the extension (i.e. the higher-order counterpart of the space 
$\D^{1; 1-2s}(\R^{n+1}_+)$ above) is not available yet. 

The main achievements of the present work is to fill this gap. 
We stress the fact that the understanding of the more appropriate functional framework constituted the most important and hardest step in our investigation, since our point of view is in fact purely analytic.

Let us now describe our approach; we refer to the next section for preliminaries and notation.

\medskip

In order to face certain severe technical obstructions (which do not occur when $s<1$) we prefer to extend the function $\E_{s}[u]$ to $\R^{n+1}$ by symmetry. More precisely, we put
$$
\E_{s}[u](x,y) = 
\begin{cases}
(u* \P_{\!s}^y)(x)& \text{if $y\neq 0$}\\
u(x)& \text{if $y= 0$,}
\end{cases}
$$
so that $\E_{s}[u]$ is even in the $y$-variable.
Since $\mathcal D^s(\R^n)\hookrightarrow L^{\frac{2n}{n-2s}}(\R^n)$ by Sobolev embedding theorem and
$\P^y_s\in C^\infty(\R^n)\cap  L^1(\R^n)$ for any fixed $y$, then $\E_{s}[u]$ is 
well defined and measurable on $\R^{n+1}$ for any  $u\in \mathcal D^s(\R^n)$, and moreover $\E_{s}[u](\cdot,y)\in C^\infty(\R^n)\cap L^{\frac{2n}{n-2s}}(\R^n)$ provided that $y\neq 0$.

In Section \ref{S:spaces} we  introduce and study a large family of Hilbert spaces which
includes 
$$
{\mathcal D}^{1+[s];b}_{\me}(\R^{n+1})~,\qquad \|U\|_{{\mathcal D}^{1+[s];b}_{\me}}^2 =\iint\limits_{\R^{n+1}}|y|^{b}|\bnabla_{\!b}^{1+[s]}U|^2~\!dz~\!,
$$
such that any  $U\in {\mathcal D}^{1+[s];b}_{\me}(\R^{n+1})$ is \underline{even} in the  $y$
variable.
Here 
$$
b:=1-2(s-[s])~,\qquad \bnabla_{\!b}^{1+[s]} = 
\begin{cases}
\bDelta_b^{\!\!\frac{1+[s]}{2}} & \text{if $[s]$ is odd}\\
\nabla\bDelta_b^{\!\!\frac{[s]}{2}} & \text{if $[s]$ is even}
\end{cases}~,\qquad  \bDelta_b = \Delta + by ^{-1}\partial_y~\!.
$$
Notice that $b\in(-1,1)$. If $s+\frac12$ is an integer, then $b=0$ and $(-\bDelta_b)^{1+[s]}$ coincides with the standard polyharmonic operator of order 
$1+2s$, which is the
$(\frac12+s)$-th power of the  Laplacian $-\Delta$. If $s$ is not an half integer, then the operator
$-\bDelta_b=-\Delta-by^{-1}\partial_y$ has a singularity at $\{y=0\}$; nevertheless, it works smoothly on functions on $\R^{n+1}$ which are even in the $y$-variable,
see Section \ref{S:Preliminaries}.

Differently from \cite{CdMG, Case}, for instance, in our approach Hardy type inequalities play a crucial role and constitute 
the first step in our investigation.
As a consequence of the more general Theorem \ref{T:Hardy}, which might have an independent interest,
one immediately obtains the next result (see Remark \ref{R:Hardy_references} for related references).

\begin{Theorem}
\label{T:Hardy_vero}
Let $s\in(0,n/2)$ be not an integer.  If $U\in {\mathcal D}^{1+[s];b}_{\me}(\R^{n+1})$, then
$$
\iirn|y|^{b}|\bnabla_{\!b}^{1+[s]}U|^2~\!dz\ge 2^{2(1+[s])}
\frac{\Gamma\big(\frac{n-2s}{4}+1+[s]\big)^2}{\Gamma\big(\frac{n-2s}{4}\big)^2}~\iirn |y|^{b}|z|^{-2(1+[s])}|U|^2~\!dz~\!.
$$
\end{Theorem}

\smallskip

Next, it turns out that any function 
$U\in {\mathcal D}^{1+[s];b}_{\me}(\R^{n+1})$ has a trace ${\rm Tr}(U)=U_{|\{y=0\}}\in \mathcal D^s(\R^n)$
and that the trace map
$$
{\rm Tr}:{\mathcal D}^{1+[s];b}_{\me}(\R^{n+1})\to \mathcal D^s(\R^n)
$$
is continuous. This  important
and difficult result is proved in Subsection \ref{SS:traces}. (see also \cite{FF} for $s\in(1,2)$).  
Moreover, the following facts hold,
\begin{itemize}
\item[$T1)$] the norm operator of the trace map ${\rm Tr}$  is given by $\|{\rm Tr}\|^2_{\mathcal L}=\frac{1}{2d_s}$, where
\begin{equation}
\label{eq:ds}
d_s= \frac{[s]!}{\Gamma(s)} \,  2^{1-2(s-[s])}\Gamma(1- (s-[s]))~\!; 
\end{equation}
\item[$T2)$] its 
adjoint $\textrm{\rm Tr}^*:\mathcal D^s(\R^n)\to \mathcal D^{1+[s];b}_\me(\R^{n+1})$ is proportional to the extension operator $\E_s$, precisely
$$
\textrm{\rm Tr}^*= \frac{1}{2d_s}\E_s~\!.
$$
\end{itemize}
Properties $T1)$ and $T2)$ readily follows from the next theorem, in which we summarize our results about $s$-polyharmonic extensions.

\begin{Theorem}
\label{T:main}
Assume that $s\in(0,n/2)$ is not an integer, and let $u\in {\mathcal D}^s(\R^n)$.
Then 
\begin{itemize}
\item[$i)$] $\E_s[u]\in {\mathcal D}^{1+[s];b}_{\me}(\R^{n+1})$ and $\text{\rm Tr}(\E_s[u])=u$;
\item[$ii)$] The function $\E_s[u]$ is the unique solution to the convex minimization problem
\begin{equation}
\label{eq:min1}
\inf_{V\in \mathcal D^{1+[s];b}_\me(\R^{n+1})}\frac12\iirn|y|^b|\bnabla_{\!b}^{1+[s]} V|^2~\!dz-
2d_s\langle \Ds u,\textrm{\rm Tr}(V)\rangle.
\end{equation}
In particular, 
\begin{equation}
\label{eq:aggiunto}
\iirn|y|^b\bnabla_{\!b}^{1+[s]} \E_s[u]~\bnabla_{\!b}^{1+[s]}V~\!dz=2d_s\langle\Ds u,\textrm{\rm Tr}(V)\rangle,\quad 
\text{for any $V\in \mathcal D^{1+[s];b}_\me(\R^{n+1})$,}
\end{equation}
and therefore
\begin{gather}
\label{eq:uguale_norme}
\iint\limits_{\R^{n+1}}|y|^b|\bnabla_{\!b}^{1+[s]}\E_s[u]|^2~\!dz=2d_s\irn|\Dshalf u|^2~\!dx~\!,\\
\nonumber
-\div(|y|^b\nabla(-\bDelta_b)^{[s]} \E_s[u])~=0\qquad\text{on~ $\R^{n+1}_+$;}
\end{gather}
\item[$iii)$] The function $\E_s[u]$ is the unique solution to the convex minimization problem
\begin{equation}
\label{eq:min2}
\displaystyle \inf_{V \in {\mathcal D}^{1+[s];b}_{\me}(\R^{n+1}) \atop \textrm{\rm Tr}(V)=u} \iint\limits_{\R^{n+1}}|y|^b|\bnabla_{\!b}^{1+[s]}V|^2~\!dz~\!;
\end{equation}
\item[$iv)$] 
 It holds that
\begin{equation}
\label{eq:uguale_operatori}
-\lim_{y \to 0}|y|^{b-1}y~\!\frac{\partial}{\partial y} (-\bDelta_b)^{[s]}\E_s[u] =  d_s\,(-\Delta)^s u
\end{equation}
in the dual space $\mathcal D^{-s}(\R^n)$. If in addition $s>1$, then for any integer $m=1,\dots,[s]$ we have that
\begin{gather}
\label{eq:andata}
(-\bDelta_b)^{m}\E_s[u] = -2(1+[s]-m)~y^{-1}\frac{\partial}{\partial y}(-\bDelta_b)^{m-1}\E_s[u]\quad\text{on $\{y\neq 0\}$}
\\
\label{eq:limits}
\lim_{y \to 0} (-\bDelta_b)^{m}\E_s[u] =
 \frac{d_s}{d_{s-m}}(-\Delta)^{m} u\\
{\lim\limits_{y\to 0} y^{-1}\frac{\partial^{2m-1}}{\partial y^{2{m}-1}} \E_s[u]} =
\lim\limits_{y\to 0}
{\frac{\partial^{2m}}{\partial y^{2m}}\E_s[u] }= \kappa_{s,{m}} (-\Delta)^{m} u
\label{eq:limits2}
\end{gather}
for some explicit constant $\kappa_{s,m}$ (see Lemma 
\ref{L:extreg1}). The above limits are taken in the sense of traces if $s>2m$ and in the dual space $\mathcal D^{s-2m}(\R^n)$ if $s<2m$;
\item[$v)$] If $u\in C^\infty_c(\R^n)$, then $\E_s[u]\in  C^{2[s],\sigma}(\R^{n+1})$
for any $\sigma\in(0,1)$, and the limits in  $iv)$ hold in the uniform topology of $\R^n$.
\end{itemize}
\end{Theorem}

\medskip

\noindent
Our proof of Theorem \ref{T:main} in Section \ref{S:proof} needs some preliminary work.
Besides the already mentioned Hardy type inequalities and trace theorems, a careful investigation of the \underline{family} of Poisson kernels $\{\P^y_\alpha\}_{\alpha>0}$ 
and relative extension operators $\{\E_\alpha\}_{\alpha>0}$, together
with the relations within them, are needed as well (see Section \ref{S:genext}).

We stress the fact that Theorem \ref{T:main} is essentially known if $s\in(0,1)$, see  \cite{CS, FrLe, FLS}; some of its conclusions
can be recovered by using the results in  \cite{Case}.

Few  comments are in order. 
One can reformulate $ii)$ in Theorem \ref{T:main} by saying that
$\E_s[u]$ is a weak solution to
\begin{equation}\label{eq:0_extequation}
\begin{cases}
-\div(|y|^b\nabla(-\bDelta_b)^{[s]} U) = 2d_s~\!\delta_{\{y=0\}}\Ds u & \text{on }\R^{n+1}\\
U(\cdot,0)=u~\!,
\end{cases}
\end{equation}
see Subsection \ref{SS:distribuzionale} for definitions and details.

Next, let $u\in C^\infty_c(\R^n)$. Thanks to $v)$ in Theorem \ref{T:main}, we can write the Taylor expansion formula
$$
\E_s[u](\cdot,y)=\sum_{m=0}^{[s]} \frac{\kappa_{s,m}}{(2m)!}~\!y^{2m}~\!(-\Delta)^m u+o(y^{2[s]})\quad \text{
uniformly on $\R^n$, as $y\to 0$.}
$$

Finally, 
Theorem \ref{T:main} gives  informations about polyharmonic extensions to the upper half space of any 
sufficiently regular function $u$ on $\R^n$.
In fact, if $k\ge 1$, by choosing $s=k-\frac12$ one has  $b=0$ and $(-\bDelta_b)^{1+[s]}=(-\Delta)^k$. Hence $U:=\E_{k-\frac12}[u]$ is the unique weak (i.e. in a suitable energy space) solution to
$$
(-\Delta)^kU=0\qquad \text{on $\R^{n+1}_+$,}\qquad U(\cdot,0)=u.
$$

\medskip
To introduce our last main result we recall the Hardy inequality by Herbst \cite{He},
$$
\irn|\Dshalf u|^2~\!dx\ge 2^{2s}~\!\frac{\Gamma\big(\frac{n+2s}{4}\big)^2}{\Gamma\big(\frac{n-2s}{4}\big)^2}~\irn|x|^{-2s}|u|^2~\!dx~
\quad\text{for any $u\in \mathcal D^s(\R^n)$.}
$$
In the next Theorem we generalize \cite[Lemma 2.1]{MN1},  \cite[Theorem 1]{MN2}
to higher orders, and  give a positive answer to a question raised in \cite[Remark 2.2]{MN1} for $n=1, s\in(0,\frac12)$.
\begin{Theorem}
\label{T:mia}
Let $s \in (0, n/2)$. Then 
$$
\iint\limits_{\R^{n+1}}|y|^{1-2(s-[s])}|z|^{-2(1+[s])}|\E_s[u]|^2~\!dz \le  \gamma \irn|x|^{-2s}|u|^2~\!dx~\quad
\text{for any $u\in \mathcal D^s(\R^n)$,}
$$
where the positive constant $\gamma$ does not depend on $u$. 
\end{Theorem}
Theorem \ref{T:mia} is an immediate consequence of Lemma \ref{L:L2_continuita}, which provides similar estimates for the extension operators $\E_\alpha$ for any $\alpha>0$.

Differently from the arguments in \cite{MN1, MN2}, the proof of Theorem \ref{T:mia} relies on the characterization of Muckenhoupt weights via the Hardy-Littlewood maximal operator. We believe that it can be further generalized in order to consider additional parameters, in the spirit of \cite{MN2}; nevertheless, this is beyond the aim of the present work.

\section{Notation and preliminaries}
\label{S:Preliminaries}

 We start by listing some notations used throughout the paper.

 $\bullet$ $\R^{n+1} \equiv \R^n \times \R =  \{ z=(x,y) ~|~ x \in \R^n, ~y \in \R\}~,~~
 \R^{ n+1}_+ =\R^n \times (0,\infty)$;
 
 $\bullet$ If $\zeta \in \R^d$ and $r>0$, then $B_r(\zeta)$ is the ball of radius $r$ about $\zeta$ in $\R^d$;

$\bullet$ $dz = dx dy$ is the volume element in $\R^{n+1}$; 

$\bullet$ Let $k\ge 0$ be an integer or $k=\infty$. We put 
$$
C^k_{\me}(\R^{n+1})=\{\ U\in C^k(\R^{n+1})~|~ U(x,\cdot)~~\text{is even}~\}
~,\quad C_{c;\me}^k(\R^{n+1}) = C_{\me}^k(\R^{n+1})\cap C^k_c(\R^{n+1})
$$
and regard at $C^k_{c;\me}(\R^{n+1})$ as a subspace of $(C^k_{c}(\R^{n+1}), \|\cdot\|_{C^k_c})$;

$\bullet$ We endow $L^p(\R^d)$ with the standard norm $\|\cdot\|_p$~\!. If
 $\omega\ge 0$ is a measurable function on $\R^d$, then the weighted space $L^2(\R^d; \omega d \zeta)$ 
inherits an Hilbertian structure with respect to the  norm $\|\omega^{1/2}u\|_2$;

$\bullet$ The $0$-th power of any differential operator is the identity; 

$\bullet$ $\displaystyle \partial ^j_y := \frac{\partial^j}{\partial y^j}$ for any integer $j \geq 0$; 

$\bullet$ By $c$ we denote generic constants, whose value may change from line to line.

\paragraph{Muckenhoupt weights.}
We denote by $A_2(\R^d)$ the class of Muckenhoupt weights on $\R^d$, which are nonnegative functions 
$\omega\in L^1_{\rm loc}(\R^d)$ such that
\begin{equation}
\label{eq:comega}
c_\omega:=\sup_{\zeta \in \R^d,\ r>0}\Big(\frac{1}{|B_r(\zeta)|} \int\limits_{B_r(\zeta)} \omega\,d\xi \Big) \Big(\frac{1}{|B_r(\zeta)|} \int\limits_{B_r(\zeta)}\omega^{-1} \,d\xi \Big) < \infty\,.
\end{equation}
It is known that $\omega\in A_2(\R^d)$ if and only if the Hardy--Littlewood maximal operator
\begin{equation}\label{eq:HLop}
M_d[u](\zeta) = \sup_{r>0}\frac{1}{|B_r(\zeta)|}\int\limits_{B_r(\zeta)} |u(\xi)|\,d\xi\,, 
\end{equation}
is bounded $L^2(\R^d;\omega \,d\zeta)\to L^2(\R^d;\omega \,d\zeta)$. 

\paragraph{Fractional Laplacian and Fourier transform.}

The fractional Laplacian $\Ds$ of a rapidly decaying function $u$ on $\R^n$ is  defined 
via the Fourier transform by 
$$
\widehat{\Ds u} = |\xi|^{2s}\widehat{u}~\!, \qquad
\widehat{u}(\xi)= (2\pi)^{-\frac{n}{2}}\int\limits_{\R^n} e^{-i~\!\!\xi\cdot x}u(x)~\!dx\,.
$$
Let $n>2s>0$. Thanks to the Hardy inequality \cite{He}, the space
$$
\mathcal D^s(\R^n)=\big\{u\in L^2(\R^n; |x|^{-2s}dx)~|~\Dshalf u \in L^2(\R^n)~\!\big\}
$$
naturally inherits a Hilbertian structure from the scalar product
$$
(u,v)=\irn\Dshalf u\Dshalf v~\!dx=\irn|\xi|^{2s}\,\widehat{u}\,\overline{\widehat{v}}~\!d\xi\,.
$$
It is well known that smooth, compactly supported functions are dense in $\mathcal D^s(\R^n)$.

\paragraph{Weighted polyharmonic operators.} 
Given any integer $j\ge 1$ and any $b \in (-1,1)$, we formally introduce the following differential operators 
for functions on $\R^{n+1}$, 
$$
\bDelta_b = \Delta + by ^{-1}\partial_y,\quad 
\bnabla_{\!b}^{j} = 
\begin{cases}
\bDelta_b^{\!\!\frac{j}{2}} & \text{if $j$ is even}\\
\nabla\bDelta_b^{\!\!\frac{j-1}{2}} & \text{if $j$ is odd}. 
\end{cases}
$$

\medskip
If $U\in C^2_{\me}(\R^{n+1})$,  then 
$y^{-1}\partial_yU(x,y)=\partial^2_yU(x,0)+o(1)$ as $y\to 0$, uniformly for $x$ on compact sets of $\R^n$.
Thus, $-\bDelta_b U \in C^{0}_{\me}(\R^{n+1})$. More generally, 
\begin{equation}\label{eq:infreg}
\begin{cases}
\bnabla_{\!b}^{j} U \in C^{k-{j}}_\me(\R^{n+1})&\text{if ${j}$ is even,}\\
\bnabla_{\!b}^{j} U \in C^{k-{j}}(\R^{n+1})^{n+1}&\text{if ${j}$ is odd,}
\end{cases}
\quad \text{for any ~$U \in C^{k}_\me(\R^{n+1})\,$, $\,{j}=0,\dots,k$\,.}
\end{equation}
If in addition $U\in C^2_{c;\me}(\R^{n+1})$ has compact support, then $\|\bDelta_bU\|_\infty\le c(b,U)\|U\|_{C^2_c}$ where $c(b,U)$ depends on the support of $U$.
Using induction, it is easy to infer that 
\begin{equation}
\label{eq:that}
\|\bnabla_{\!b}^k U\|_\infty\le c(b,U)\|U\|_{C^k_c}~,\quad \text{for any ~$U \in C^{k}_{c;\me}(\R^{n+1})$}\,.
\end{equation}

Let $j, h\ge 1$ be odd. With some abuse of notation, for $\f \in C^{j}_{\me}(\R^{n+1})$, $\psi\in C^{h}_{\me}(\R^{n+1})$, we  put 
$$
\bnabla_{\!b}^{j}\f\bnabla_{\!b}^{h}\psi := (\nabla\bDelta_b^{\frac{{j}-1}{2}}\!\f)\cdot(\nabla\bDelta_b^{\frac{{h}-1}{2}}).
$$

We point out a useful integration by parts formula.
\begin{Lemma}
\label{L:Step1}
Let $k\ge 2$, $b\in (-1,1)$, $W\in C^{2(k-1)}_\me(\R^{n+1})$. Then
\begin{equation}
\label{eq:basta}
\iint\limits_{\R^{n+1}} |y|^b(-\bDelta_b)^{{k-1}}W~\!(-\bDelta_b) V~\!dz= \iint\limits_{\R^{n+1}} |y|^b\bnabla_{\!b}^{{k}}W~
\bnabla_{\!b}^{{k}} V~\!dz\quad \text{for any $V\in C^\infty_{c;\me}(\R^{n+1})$.}
\end{equation}
\end{Lemma}

\proof
In this proof we neglect to write the volume integration form $dz$. 

Fix $V\in C^\infty_{c;\me}(\R^{n+1})$. Notice that $(-\bDelta_b)^{{k-1}}W$ and $\bnabla_{\!b}^{{k}}W$ are continuous functions on $\R^{n+1}$. Since $\bDelta_b V$ and $\bnabla_{\!b}^{{k}} V$ are smooth and compactly supported, 
the integrals in (\ref{eq:basta}) are well defined and finite.

If ${k}=2m+1$ is odd, then $W\in C^{4m}_\me(\R^{n+1})$. We have to prove that 
\begin{equation}
 \label{eq:g_odd}
\iint\limits_{\R^{n+1}} |y|^b\bDelta_b^{2m}W~\!\bDelta_b V=-
\iint\limits_{\R^{n+1}} |y|^b \nabla(\bDelta_b^{m}{W})\cdot \nabla(\bDelta_b^{m}{V})\,.
 \end{equation}
 If $W\in C^4_\me(\R^{n+1})$, then $\bDelta_b{W}\in C^2_\me(\R^{n+1})$ by (\ref{eq:infreg}). Since
 $\bDelta_b{V}\in C^\infty_{c;\me}(\R^{n+1})$,
 we can use integration
by parts to obtain
$$
\iint\limits_{\R^{n+1}} |y|^b\bDelta^2_b{W}~\!\bDelta_b{V}=
\iint\limits_{\R^{n+1}} \div(|y|^b\nabla {\bDelta_bW})(\bDelta_b{V})
=- \iint\limits_{\R^{n+1}} |y|^b\nabla {\bDelta_bW}\cdot\nabla\bDelta_b{V}~\!.
$$
We proved that  
\begin{equation}
\label{eq:even_2}
\iint\limits_{\R^{n+1}}|y|^b\bDelta^2_b{W}~\!\bDelta_b{V}=-
\iint\limits_{\R^{n+1}}|y|^b\nabla {\bDelta_bW}\cdot\nabla\bDelta_b{V}
\quad \text{for any $W\in C^4_\me(\R^{n+1}), V\in C^\infty_{c;\me}(\R^{n+1})$}\,.
\end{equation}
In particular, if $m=1$ then (\ref{eq:g_odd}) follows. 
 
Assume now that (\ref{eq:g_odd}) holds
for some integer $m$. If $W\in C^{4(m+1)}_\me(\R^{n+1})$ we have
$$
\begin{aligned}
\iint\limits_{\R^{n+1}} &|y|^b\bDelta_b^{2(m+1)}W~\!\bDelta_bV=
\iint\limits_{\R^{n+1}} |y|^b\bDelta_b^{2m}(\bDelta^2_b W)~\!\bDelta_bV
=-\iint\limits_{\R^{n+1}} |y|^b\nabla\big(\bDelta_b^{m+2}W\big)\cdot\nabla(\bDelta_b^m V)\\
&=
\iint\limits_{\R^{n+1}} |y|^b \bDelta_b^{m+2}W~\!\bDelta_b^{m+1} V=
\iint\limits_{\R^{n+1}} |y|^b \bDelta_b^{2}\big(\bDelta_b^mW\big)~\!\bDelta_b^{m+1} V
=- \iint\limits_{\R^{n+1}} |y|^b\nabla\big(\bDelta_b^{m+1}W\big)\cdot\nabla(\bDelta_b^{m+1} V)
\end{aligned}
$$
by (\ref{eq:even_2}), with $\bDelta_b^{m}W$ instead of $W$ and $\bDelta_b^{m} V$ instead of $V$. The "odd" case is complete.

\medskip

We now deal with the case $k=2m$, $m\ge 1$. We have to prove that 
\begin{equation}
\label{eq:g_even}
\iint\limits_{\R^{n+1}} |y|^b\bDelta_b^{2m-1}W~\!\bDelta_b V=
\iint\limits_{\R^{n+1}} |y|^b \bDelta_b^{m}{W} \, \bDelta_b^{m}{V}
\quad \text{for any $W\in C^{2(2m-1)}_\me(\R^{n+1}), V\in C^\infty_{c;\me}(\R^{n+1})$.}
 \end{equation}
The case $m=1$ is trivial. Assume  that (\ref{eq:g_even}) holds for some integer $m\ge 1$
and let $W\in C^{2(2m+1)}_{c;\me}(\R^{n+1})$. Since $\bDelta_b^{2(m+1)-1} W=
\bDelta_b^{2m-1}\big(\bDelta_b^2W)$, using (\ref{eq:g_even}) and then (\ref{eq:even_2}) we obtain
$$
\begin{aligned}
\iint\limits_{\R^{n+1}} |y|^b\bDelta_b^{2(m+1)-1}W~\!\bDelta_b V&=
\iint\limits_{\R^{n+1}} |y|^b\bDelta_b^{2}\big(\bDelta_b^mW\big)~\!\bDelta_b^m V=-
\iint\limits_{\R^{n+1}}|y|^b\nabla {\bDelta^{m+1}_bW}\cdot\nabla\bDelta_b^{m}{V}\\
&=\iint\limits_{\R^{n+1}}|y|^b {\bDelta^{m+1}_bW}\bDelta_b^{m+1}{V}~\!,
\end{aligned}
$$
which concludes the proof.
\QED

\section{A class of homogeneous function spaces}
\label{S:spaces}

We need to define a large class of spaces ${\mathcal D}^{k;a,b}_{\me}(\R^{n+1})$ depending on the 
extra parameter $a\ge 0$. 

In this section  $k\ge 0$ is  integer and the fixed exponents $a,b$ satisfy
\begin{equation}\label{eq:abcond}
-1<b<1, \quad\quad \quad  
0\le a < \frac{n+1+b}{2}-k\,.
\end{equation}

\begin{Remark}
\label{R:M}
Under the above assumptions, the weights
\[
\omega(z) = |y|^b|z|^{-2(a+k-j)}=\frac{|y|^b}{(|x|^2+y^2)^{a+k-j}}~,\qquad j=0,\ldots,k\,,
\]
belong to the Muckenhoupt class $A_2(\R^{n+1})$. In fact,  
the supremum $c_\omega$ in (\ref{eq:comega}) is estimated by
\[
c_{\omega} \leq c_n\!\!\!\sup_{(x,y)\in \R^{n+1}\!,\atop r>0} \!\!r^{-2(n+1)}
\Big( \int\limits_{y-r}^{y+r}|\tau|^b~\!d\tau\!\!\!\int\limits_{B_r(x)}\! \frac{1}{(|\zeta|^2+\tau^2)^{a+k-j}}~\!d\zeta\Big)
\Big( \int\limits_{y-r}^{y+r}|\tau|^{-b}~\!d\tau\!\!\!\int\limits_{B_r(x)}\!(|\zeta|^2+\tau^2)^{a+k-j} {d\zeta}\Big)
<\infty~\!.
\]
\end{Remark}

Recall that for any $U\in C^\infty_{c;\me}(\R^{n+1})$ we have $|\bnabla_{\!b}^j U|\in C^0_c(\R^{n+1})$
for any $j=0,\dots,k$ by \eqref{eq:infreg}. Thus 
\begin{equation}
\label{eq:nn}
\nn U\nn_{k;a,b}^2 := \sum_{j=0}^{k}\ \iint\limits_{\R^{n+1}}|y|^b|z|^{-2(a+k-j)}|\bnabla_{\!b}^j U|^2dz<\infty
~\!.
\end{equation}

\begin{Definition}
The space
$$
{\mathcal D}^{k;a,b}_{\me}(\R^{n+1})
$$
is the completion of $C^\infty_{c;\me}(\R^{n+1})$ in $L^2(\R^{n+1};|y|^b|z|^{-2(a+k)}~\!dz)$
with respect to the norm $\nn \cdot\nn_{k;a,b}$.
\end{Definition}

One can show via standard arguments that ${\mathcal D}^{k;a,b}_{\me}(\R^{n+1})$ is in fact 
an Hilbert space. 

\begin{Remark}
\label{R:inclusioni}
Let $U\in \mathcal D^{k;a,b}_\me(\R^{n+1})$. Then $\bnabla^j_bU$ is defined by density for any $j=1,\dots,k$, and
$|\bnabla^j_bU|\in L^2(\R^{n+1};|y|^b|z|^{-2(a+k-j)}dz)$. Moreover,
$$
\bDelta_b^mU\in \mathcal D_\me^{k-2m;a,b}(\R^{n+1})\quad\text{for any integer $1\leq m<\frac{k}{2}$.}
$$
\end{Remark}

\begin{Remark}
\label{R:ck_dense}
It turns out that $C^{k}_{c;\me}(\R^{n+1})\subset \mathcal D_\me^{k;a,b}(\R^{n+1})$. For the proof, take 
$U \in C^{k}_{c;\me}(\R^{n+1})$ and a sequence $(\rho_{h})_{h \in \N} \subset C^{\infty}_{c;\me}(\R^{n+1})$  of radially symmetric mollifiers. Since 
$$
\|\bnabla_{\!b}^k(U*\rho_h)-\bnabla_{\!b}^k U\|_\infty\le \|U*\rho_h-U\|_{C^k_c},
$$
by (\ref{eq:that}), then $U*\rho_h\to U$ in $\mathcal D^{k;b}_{\me}(\R^{n+1})$.
\end{Remark}

Actually, we will endow ${\mathcal D}^{k;a,b}_{\me}(\R^{n+1})$ with the
more natural norm
\[
\|U\|_{k;a,b}^2=\iint\limits_{\R^{n+1}}|y|^b|z|^{-2a}|\bnabla^k_{\!b} U|^2~\!dz~\!,
\]
which  turns out to be equivalent to $\nn \cdot \nn_{k;a,b}$ thanks to the 
Hardy-type inequalities in the next result.

\begin{Theorem}[Hardy inequalities]
\label{T:Hardy}
Let $k,a,b$ as in (\ref{eq:abcond}). Then
$\|\cdot \|_{k;a,b}$ is an equivalent Hilbertian norm in $\mathcal D^{k;a,b}_\me(\R^{n+1})$. Moreover,  
\begin{equation}\label{eq:hardyN0}
\iint\limits_{\R^{n+1}}|y|^b |z|^{-2a}|\bnabla_{\!b}^{k} U|^2~\!dz\ge 
\mathcal H_{k,a,b}^2\iint\limits_{\R^{n+1}}|y|^b|z|^{-2(a+k)}|U|^2~\!dz
\quad\text{for any $U\in \mathcal D^{k;b}_\me(\R^{n+1})$,}
\end{equation}
where
$$
\mathcal H_{k,a,b}=2^k~\!\dfrac{\Gamma\big(\frac{n+1+b}{4}+\frac{k}{2}-\big[\frac{k}{2}\big]-\frac{a}{2}\big)}
{\Gamma\big(\frac{n+1+b}{4}+\frac{k}{2}-\big[\frac{k}{2}\big]+\frac{a}{2}\big)}~\!
\dfrac{\Gamma\big(\frac{n+1+b}{4}+\frac{k}{2}+\frac{a}{2}\big)}{\Gamma\big(\frac{n+1+b}{4}-\frac{k}{2}-\frac{a}{2}\big)}~\!.
$$
\end{Theorem}

\proof 
In this proof we neglect to write the volume integration form $dz$. 

Fix a nontrivial $U\in C^\infty_{c;\me}(\R^{n+1})$. We will use induction to prove (\ref{eq:hardyN0}) 
and the existence of a constant $C_{a,k}$ not depending on $U$, such that
\begin{equation}\label{eq:hardyN}
C^2_{a,k}\nn U\nn_{k;a,b}^2\le \iint\limits_{\R^{n+1}}|y|^b |z|^{-2a}|\bnabla_{\!b}^{k} U|^2\,.
\end{equation}

\paragraph{Step 1.} Let $k=1$. Then
$\displaystyle{H_1:=\mathcal H_{1,a,b}=2~\frac{\Gamma\big(\frac{n+1+b}{4}-\frac{a}{2}+\frac12\big)}{\Gamma\big(\frac{n+1+b}{4}-\frac{a}{2}-\frac12\big)}=
\frac{n+1+b}{2}-(a+1)>0.
}$
We have to prove that
\begin{equation}\label{eq:hardy1}
{H^2_1} \iint\limits_{\R^{n+1}}|y|^b|z|^{-2(a+1)}|U|^2 \leq \iint\limits_{\R^{n+1}}|y|^b|z|^{-2a}|\nabla U|^2.
\end{equation}
We can assume that $a>0$; the case
$a=0$ is easily recovered by taking the limit as $a\searrow 0$.

We have $|y|^b|z|^{-2a}, |y|^b|z|^{-2(a+1)} \in L^1_{loc}(\R^{n+1})$ and
\begin{equation}
\label{eq:za1}
-\div(|y|^b\nabla|z|^{-2a}) =  4a{H_1} |y|^b|z|^{-2(a+1)}\quad \text{on $\{y\neq 0\}$}\,.
\end{equation}
In fact, (\ref{eq:za1}) holds true in the distributional sense, and thanks to a standard
approximation argument we can use integration by parts and  H\"older inequality to obtain
$$
4a{H_1}
\!\!\!\iint\limits_{\R^{n+1}}|y|^b|z|^{-2(a+1)}|U|^2\!=
\!\!\!\iint\limits_{\R^{n+1}}|y|^b\nabla |z|^{-2a}\cdot \nabla U^2
\leq 
4a \Big[\!\!\!\iint\limits_{\R^{n+1}}|y|^b|z|^{-2(a+1)}|U|^2\Big]^{\!\!\frac{1}{2}}\Big[\!\!\!\iint\limits_{\R^{n+1}}|y|^b|z|^{-2a}|\nabla U|^2\Big]^{\!\!\frac{1}{2}}.
$$
The conclusion readily follows, as (\ref{eq:hardyN0}) and (\ref{eq:hardyN}) are equivalent in this case.

\paragraph{Step 2.} Let $k=2$. Then
$$
\mathcal H_{2,a,b}=2^2~\!
\dfrac{\Gamma\big(\frac{n+1+b}{4}-\frac{a}{2}\big)}{\Gamma\big(\frac{n+1+b}{4}+\frac{a}{2}\big)}~\!
\frac{\Gamma\big(\frac{n+1+b}{4}+\frac{a}{2}+1\big)}{\Gamma\big(\frac{n+1+b}{4}-\frac{a}{2}-1\big)}=
(H_1-1)\big(\frac{n+1+b}{2}+a\big).
$$
Now the starting point is the equality 
(\ref{eq:za1}) with $a$ replaced by $a+1$, that is,
$$
-\div(|y|^b\nabla|z|^{-2(a+1)}) =  4(a+1)(H_1-1) |y|^b|z|^{-2(a+2)}\quad \text{on $\{y\neq 0\}$}\, .
$$
Since $U(x,\cdot)$ is even, we can
integrate by parts two times to obtain
\[
\begin{aligned}
-4(a+1)(H_1-1)& \!\!\iint\limits_{\R^{n+1}}|y|^b|z|^{-2(a+2)}|U|^2=
\!\!\iint\limits_{\R^{n+1}}\div(|y|^b\nabla |z|^{-2(a+1)})|U|^2
=\!\!\iint\limits_{\R^{n+1}}|z|^{-2(a+1)}\div(|y|^b\nabla U^2)\\
&=
\!\!\iint\limits_{\R^{n+1}}|y|^b|z|^{-2(a+1)}\bDelta_bU^2
= 2\!\!\iint\limits_{\R^{n+1}}|y|^b|z|^{-2(a+1)}U\bDelta_bU+
2\!\!\iint\limits_{\R^{n+1}}|y|^b|z|^{-2(a+1)}|\nabla U|^2~\!.
\end{aligned}
\]
Thanks to  H\"older inequality we infer that 
\begin{equation}
\label{eq:h2}
\begin{aligned}
\iint\limits_{\R^{n+1}}|y|^b|z|^{-2(a+1)}|\nabla U|^2&+
2(a+1)(H_1-1) \iint\limits_{\R^{n+1}}|y|^b|z|^{-2(a+2)}|U|^2\\
&\le \Big(\iint\limits_{\R^{n+1}}|y|^b|z|^{-2(a+2)}| U|^2\Big)^\frac12
\Big(\iint\limits_{\R^{n+1}}|y|^b|z|^{-2a}| \bDelta_b U|^2\Big)^\frac12\,.
\end{aligned}
\end{equation}
Now we use (\ref{eq:hardy1}) with $a$ replaced by $a+1$ to estimate
\begin{equation}
\label{eq:h2bis}
\iint\limits_{\R^{n+1}}|y|^b|z|^{-2(a+1)}|\nabla U|^2\ge 
(H_1-1)^2 \iint\limits_{\R^{n+1}}|y|^b|z|^{-2(a+2)}|U|^2
\end{equation}
which, together with (\ref{eq:h2}), gives
$$
\mathcal H_{2,a,b}^2\iint\limits_{\R^{n+1}}|y|^b|z|^{-2(a+2)}|U|^2
\le \iint\limits_{\R^{n+1}}|y|^b|z|^{-2a}| \bDelta_b U|^2
$$
as $(H_1-1)^2+2(a+1)(H_1-1)=\mathcal H_{2,a,b}$. 
To conclude Step 2 notice that (\ref{eq:h2}) and (\ref{eq:h2bis}) trivially imply
$$
(H_1-1)^2\iint\limits_{\R^{n+1}}|y|^b|z|^{-2(a+1)}| \nabla U|^2 \le
\iint\limits_{\R^{n+1}}|y|^b|z|^{-2a}| \bDelta_b U|^2\,.
$$

\paragraph{Step 3.} It remains to consider the case $k\ge 3$. If $k=2m+1$ is odd, then
$$
\begin{aligned}
\iint\limits_{\R^{n+1}}|y|^b|z|^{-2a}|\bnabla^k_{\!b} U|^2&=
\iint\limits_{\R^{n+1}}|y|^b|z|^{-2a}|\nabla(\bDelta_b^m U)|^2\\
&\ge {H^2_1}\iint\limits_{\R^{n+1}}|y|^b|z|^{-2(a+1)}|\bDelta_b^m U|^2
= {H^2_1}\iint\limits_{\R^{n+1}}|y|^b|z|^{-2(a+1)}|\bnabla_{\!b}^{k-1} U|^2.
\end{aligned}
$$
If $k=2m\ge 4$ is even we write two chain of inequalities,
$$
\begin{aligned}
\iint\limits_{\R^{n+1}}|y|^b|z|^{-2a}|\bnabla^k_{\!b} U|^2&=
\iint\limits_{\R^{n+1}}|y|^b|z|^{-2a}|\bDelta_b^m U|^2\\
&\ge {\mathcal H_{2,a,b}^2}\iint\limits_{\R^{n+1}}|y|^b|z|^{-2(a+2)}|\bDelta_b^{m-1} U|^2
= {\mathcal H_{2,a,b}^2}\iint\limits_{\R^{n+1}}|y|^b|z|^{-2(a+2)}|\bnabla_{\!b}^{k-2} U|^2,
\end{aligned}
$$
$$
\begin{aligned}
\iint\limits_{\R^{n+1}}|y|^b|z|^{-2a}&|\bnabla^k_{\!b} U|^2=
\iint\limits_{\R^{n+1}}|y|^b|z|^{-2a}|\bDelta_b(\bDelta_b^{m-1} U)|^2\\
&\ge (H_1-1)^2\iint\limits_{\R^{n+1}}|y|^b|z|^{-2(a+1)}|\nabla(\bDelta_b^{m-1} U)|^2
= (H_1-1)^2\iint\limits_{\R^{n+1}}|y|^b|z|^{-2(a+1)}|\bnabla_{\!b}^{k-1} U|^2~\!.
\end{aligned}
$$
The above inequalities, together with the lower order cases $k=1$ and $k=2$ and induction
easily lead to the conclusion of the proof. We omit  details.
\QED

\begin{Remark}
\label{R:Hardy_references}
If $b=0$, then  (\ref{eq:hardyN0}) reduces to the
classical (weighted) Hardy ($k=1$) and Rellich ($k=2$) inequalities. We cite also
\cite[Theorem 3.3]{Mit} for more general sharp inequalities.

The case $a=0$, $k=2$, $b\in(-1,1)$ has been already discussed in \cite{FF}.

It would be of interest to prove that the explicit constant $
\mathcal H_{k,a,b}^2$ in (\ref{eq:hardyN0})
is sharp and not achieved (this is well known in case $b=0$).
\end{Remark}

\subsection{The space $\mathcal D^{k;b}_\me(\R^{n+1})$ and trace theorems}
 \label{SS:traces}
 
We fix an integer $k\ge 1$, an exponent $b\in(-1,1)$ and focus our attention on the Hilbert space
\[
{\mathcal D}_\me^{k;b}(\R^{n+1}):={\mathcal D}_\me^{k;0, b}(\R^{n+1})
\quad \text{with norm}
\quad
\|U\|^2_{{k;b}}= \iint\limits_{\R^{n+1}}|y|^b|\bnabla_{\!b}^{k} U|^2~\!dz~\!.
\]
We start by pointing out an immediate consequence of Theorem \ref{T:Hardy}.
\begin{Corollary}
\label{C:Hardy}
Assume that $n+1+b>2k$. If $U\in {\mathcal D}^{k;b}_{\me}(\R^{n+1})$, then
$$
\iint\limits_{\R^{n+1}}|y|^b |\bnabla_{\!b}^{k} U|^2~\!dz\ge 
2^{2k}~\!\dfrac{\Gamma\big(\frac{n+1+b}{4}+\frac{k}{2}\big)^2}{\Gamma\big(\frac{n+1+b}{4}-\frac{k}{2}\big)^2}
\iint\limits_{\R^{n+1}}|y|^b|z|^{-2k}|U|^2~\!dz~\!.
$$
\end{Corollary}
The main result in this section is Theorem \ref{T:traces} below (compare with  \cite[Section 3]{FF} for $s\in(1,2)$).

\begin{Theorem}
\label{T:traces}
If $n+1+b>2k$ then the trace map $U \mapsto {\rm Tr}(U):=U(x,0)$, $U \in C^{k}_{c;\me}(\R^{n+1})$
can be uniquely extended to a continuous operator
$$
{\rm Tr}:{\mathcal D}_\me^{k;b}(\R^{n+1}) \mapsto \mathcal \D^{k-\frac{1+b}{2}}(\R^n)\, .
$$
\end{Theorem}

\proof 
Fix $U\in C^{k}_{c;\me}(\R^{n+1})$. 
We will show that 
\begin{equation}\label{eq:finalTrace}
\iint\limits_{\R^{n+1}}|y|^b|\bnabla^{k}_{\!b} U|^2~\!dz \geq c_{b}
\irn |\Dshalf U(x,0)|^2~\!dx~,\qquad s:=k-\frac{1+b}{2}~\!,
\end{equation}
where
\begin{equation}\label{eq:step2}
c_b := \inf_{\Phi \in C^1_\me(\R) \atop \Phi(0)=1} \int\limits_{-\infty}^\infty |t|^b(|\Phi'|^2 + |\Phi|^2)dt~\!.
\end{equation}
Let us first prove that $c_b>0$, which, together with (\ref{eq:finalTrace}), concludes the proof.

\medskip

By contradiction, let $c_b = 0$. For every $\e \in (0,1)$ find $\Phi_\e \in  C^1_\me(\R)$ such that $\Phi_\e(0)=1$ and
\[
\int\limits_0^\infty t^b(|\Phi'_\e|^2 + |\Phi_\e|^2)dt < \e<1\,.
\]
Let $\delta_b>0$ be defined by  $4\delta_b^{1-b}= 1-b$.
For any $t\in(0,\delta_b)$ we have
$$
|\Phi_\e(t)|\ge 1-\int\limits_0^t|\Phi_\e'|~\!d{\tau}
\ge 1-
\Big(\int_0^{\delta_b} |{\tau}|^{-b} d{\tau}\Big)^{\frac{1}{2}}\ \Big(\int_0^\infty |{\tau}|^{b}|\Phi_\e'|^2 d{\tau}\Big)^{\frac{1}{2}}\ge 1-\sqrt{\frac{\delta_b^{1-b}}{1-b}}
=\frac12~\!.
$$
Therefore
$\displaystyle{\e > \int_0^{\delta_b}t^b|\Phi_\e|^2~\!dt\ge \frac14 \int_0^{\delta_b}t^b~\!dt=
\frac{1}{4(1+b)} \delta_b^{1+b}}$,
which is a contradiction for $\e$ small enough.

\medskip 

In order to  prove that (\ref{eq:finalTrace}) holds true
we  use the Fourier transform $\widehat{U}(\cdot, y)$ 
of $U(\cdot, y)$. We have
\[
\begin{aligned}
\iint\limits_{\R^{n+1}}|y|^b|\nabla U|^2~\!dz  &=
\int\limits_{\R^n}d\xi
\!\int\limits_{-\infty}^\infty |y|^b \big(|\partial_y \widehat{U}|^2+|\xi|^2|\widehat{U}|^2\big)~\!dy
= \int\limits_{\R^n}|\xi|^{1-b}~\!d\xi\int\limits_{-\infty}^\infty |t|^b(|\phi_\xi'|^2+|\phi_\xi|^2)~\!dt\,,
\end{aligned}
\] 
where, for $\xi\neq 0$, we put
\begin{equation}
\label{eq:Uphi}
\phi_\xi(t) := \hat U (\xi, |\xi|^{-1}t)~,\quad \phi_\xi\in C^\infty_{c}(\R)\, .
\end{equation}
Since
$\displaystyle{
\int_\R |t|^b(|\phi_\xi'|^2+|\phi_\xi|^2)~\!dt \ge c_b~\!|\phi_\xi(0)|^2=c_b ~\!|\widehat{U}(\xi,0)|^2}$
by (\ref{eq:step2}), we plainly infer  
$$
\iint\limits_{\R^{n+1}}|y|^b|\nabla U|^2~\!dz \ge c_b \int\limits_{\R^n}|\xi|^{1-b}|\widehat{U}(\xi,0)|^2~\!d\xi=
c_b \int\limits_{\R^n}|(-\Delta)^{\frac{s}{2}}U(x,0)|^2~\!dx, \quad s=1-\frac{1+b}{2}\,.
$$
Thus \eqref{eq:finalTrace} holds true if $k=1$.

To handle the higher order case, for any integer $m\ge 1$ we introduce the differential operator
$L^m:C^{2m}_{c;\me}(\R)\to C^{0}_{c;\me}(\R)$, which is the $m$-th power of
\[
L\Phi = \Phi'' + bt^{-1}\Phi' - \Phi\,, \quad L: C^{2}_{c;\me}(\R) \to C^{0}_{c;\me}(\R)\,.
\]
If $\Phi \in C^2_{c;\me}(\R)$, then  $t^b\Phi(t)\Phi'(t)=o(t^{1+b})$ as $t\to 0^+$, because $|b|<1$. 
Thus we can compute
\begin{eqnarray}
\nonumber
\int\limits_0^\infty |t|^b|L\Phi|^2dt &=& \int\limits_0^\infty t^b|t^{-b}(t^b \Phi')' - \Phi|^2dt \geq \int\limits_0^\infty t^b|\Phi|^2dt-2 \int\limits_0^\infty (t^b \Phi')' \Phi~\!dt\\
&=&  \int\limits_0^\infty t^b|\Phi|^2dt+2 \int\limits_0^\infty t^b |\Phi'|^2~\!dt  \geq 
\int\limits_0^\infty t^b(|\Phi'|^2 + |\Phi|^2)~\!dt \,.
\label{eq:piove}
\end{eqnarray}
If $m\ge 2$ and $\Phi\in C^{2m}_{c;\me}(\R)$ we use (\ref{eq:piove}) to obtain
$$
\int\limits_0^\infty t^b|L^m\Phi|^2dt= \int\limits_0^\infty t^b|L(L^{m-1}\Phi)|^2dt\ge 
\int\limits_0^\infty t^b(|(L^{m-1}\Phi)'|^2 + |L^{m-1}\Phi|^2)~\!dt
\ge \int\limits_0^\infty t^b|L^{m-1}\Phi|^2~\!dt~\!,
$$
so that induction readily gives
\[
\inf_{\Phi \in C^{2m}_{c;\me}(\R) \atop \Phi(0)=1} \int\limits_{-\infty}^\infty |t|^b|L^m\Phi|^2dt 
\ge c_b\,.
\]

We are in position to complete the proof.
Using induction again one can see that 
\[
\widehat{\bDelta_b^m U}(\xi,y)=|\xi|^{2m}L^m(\phi_\xi)(|\xi|y),
\]
where the functions $U$ and $\phi$ are related by (\ref{eq:Uphi}). If $k=2m$ is even,  then
\[
\begin{aligned}
\iint\limits_{\R^{n+1}}|y|^b|\bnabla^k_{\!b} U|^2~\!dz &= 
 \int\limits_{\R^n}d\xi\int\limits_{\R}|y|^b|\widehat{\bDelta^m_{b} U}|^2 ~\!dy=
 \int\limits_{\R^n} |\xi|^{4m}d\xi \int\limits_{-\infty}^\infty |y|^b\, |L^m(\phi_\xi)(|\xi|y)|^2~\!dy \\
&= \int\limits_{\R^n}|\xi|^{2k-b-1}d\xi\int\limits_{-\infty}^\infty |t|^b |L^m(\phi_\xi)|^2~\!dt\,;
\end{aligned}
\]
if $k = 2m+1$ is odd, 
\[
\begin{aligned}
\iint\limits_{\R^{n+1}}|y|^b&|\bnabla^k_{\!b} U|^2~\!dz  = \iint\limits_{\R^{n+1}}|y|^b|\nabla \bDelta_b^{m} U|^2~\!dz  \geq \iint\limits_{\R^{n+1}}|y|^b|\nabla_{\!x}( \bDelta_b^{m} U)|^2~\!dz \\
&=\int\limits_{\R^n}|\xi|^{2(2m+1)}d\xi\!\int\limits_{-\infty}^\infty |y|^b|L^m(\phi_\xi)(|\xi|y)|^2~\!dy
= \int\limits_{\R^n}|\xi|^{2k-b-1}d\xi\int\limits_{-\infty}^\infty |t|^b |L^m(\phi_\xi)|^2~\!d t\,.
\end{aligned}
\] 
We see that, in any case,
$$
\begin{aligned}
\iint\limits_{\R^{n+1}}|y|^b|\bnabla^k_{\!b} U|^2~\!dz& 
\ge c_b \int\limits_{\R^n}|\xi|^{2k-b-1}|\phi_\xi(0)|^2~\!d\xi=
c_b \int\limits_{\R^n}|\xi|^{2k-b-1}|\widehat U(\xi,0)|^2~\!d\xi=\\&
=c_b \irn |(-\Delta)^{\!\frac{2k-b-1}{4}} U(\cdot, 0)|^2~\!dx,
\end{aligned}
$$
which concludes the proof of (\ref{eq:finalTrace}) and of the Theorem.
\QED

Theorem \ref{eq:finalTrace} together with Remark \ref{R:inclusioni} readily give the next result.

\begin{Corollary}
\label{C:traces}
Assume that $n+1+b>2k$. 
If $U\in {\mathcal D}_\me^{k;b}(\R^{n+1})$, then 
$\bDelta_b^mU\in {\mathcal D}_\me^{k-2m;b}(\R^{n+1})$ has a trace 
${\rm Tr}(\bDelta_b^mU)\in \mathcal \D^{k-2m-\frac{1+b}{2}}(\R^n)$ for any integer $m\ge 0$ such that
$k-2m\ge 1$. Moreover,
$$
\iint\limits_{\R^{n+1}}|y|^b|\bnabla_{\!b}^{k} U|^2~\!dz \geq c
\irn |\left(-\Delta\right)^{\!\frac {2(k-2m)-(1+b)}{4}}\! {\rm Tr}(\bDelta_b^mU)(x,0)|^2~\!dx
$$
where $c>0$ does not depend on $U$.
\end{Corollary}

Next we provide a further integration by parts formula, which gives some infos on
the behavior of normal derivatives on $\{y=0\}$ of functions in $\mathcal D^{k;b}_\me(\R^{n+1})$,
compare with \cite[Section 3]{FF} for $s\in(1,2)$.

\begin{Lemma}
\label{L:normal}
Let $k\ge 2$ and $U \in \D_\me^{k;b}(\R^{n+1})$. Let ${m}$ be an integer such that
$1\le m\le k/2$. Then 
\[
\lim_{y \to 0} |y|^b \partial_y \bDelta_b^{m-1} U (\cdot, y) = 0\,
\]
in a weak sense, that is, 
\begin{equation}
\label{eq:parts}
\iint\limits_{\R^{n+1}_+} y^b (\bDelta_b^{{m}} U) \f~\!dz = -\iint\limits_{\R^{n+1}_+}y^b \nabla(\bDelta_b^{m-1} U) \cdot \nabla \f~\!dz \quad \text{ for any $\f \in C^{1}_{c,\me}(\R^{n+1})$}\,.
\end{equation}
\end{Lemma}
\proof The Hardy inequalities in Theorem \ref{T:Hardy} plainly imply that the integrals in (\ref{eq:parts})
converge and depend continuously on $U\in\mathcal D^{k;b}_\me(\R^{n+1})$ for any fixed
$\f\in C_{c;\me}^1(\R^{n+1})$. Thus we can assume that $U\in C^\infty_{c;\me}(\R^{n+1})$.

Let ${m}=1$. Then the smoothness and the symmetry of $U$ in the $y$-variable 
give $y^b\partial_yU(x,y)=O(y^{1+b})$ uniformly on the 
support of $\f$, as $y\to 0^+$. Thus we can integrate by parts to get
\begin{equation}
\label{eq:parts0}
\iint\limits_{\R^{n+1}_+} y^b (\bDelta_b U) \f~\!dz = -\iint\limits_{\R^{n+1}_+}y^b \nabla U \cdot \nabla \f~\!dz.
\end{equation}
Thus (\ref{eq:parts}) holds true in this case.
If ${m}\ge 2$, it suffices to use (\ref{eq:parts0}) with $U$ replaced by $\bDelta_b^{{m-1}} U$.
 \QED

\section{Fractional Poisson kernels and  extension operators}
\label{S:genext}

In order to prove Theorem \ref{T:main} we need to study the more general class of 
Poisson kernels
\[
\P^y_\alpha (x) = c_{n,\alpha} \frac{y^{2\alpha}}{(|x|^2 + y^2)^{\frac{n+2\alpha}{2}}}\,, \quad \text{ where }\quad  c_{n,\alpha} = \frac{\Gamma\big(\frac{n+2\alpha}{2}\big)}{\pi^{\frac{n}{2}}\Gamma(\alpha)}\,,
\]
and associated extension operators 
\[
\E_{\alpha}[u](x,y) = 
\begin{cases} (u* \P_{\!\alpha}^y)(x)& \text{if $y\neq 0$}\\
u(x)& \text{if $y= 0$}
\end{cases}
\,,\quad\quad\alpha>0\,.
\]

\subsection{The Poisson kernels $\P^y_\alpha$}

Clearly, $\|\P^y_\alpha\|_{L^1(\R^n)} = 1$ for any $y\neq 0$ and $\P^y_\alpha$ is smooth on $\R^{n+1}\setminus\{y=0\}$. 
In the next simple Lemma, which will be used several times in a crucial way, we put $\Delta_x=\Delta-\partial^2_{y}$.

\begin{Lemma}
\label{L:R1}
Let $\alpha>0$ and $b\in(-1,1)$. The following identities hold on $\R^{n+1}\setminus\{y=0\}$,
\begin{itemize}
\item[$i)$] 
$\partial_y\P^y_\alpha=2\alpha y^{-1}\big(\P^y_\alpha-\P^y_{\alpha+1}\big)~,
\quad \bDelta_b \P^y_\alpha=2\alpha(b-1+2\alpha)y^{-2}\big(\P^y_\alpha-\P^y_{\alpha+1}\big)\,;$
\item[$ii)$]
$\partial_y\P^y_\alpha=\dfrac{1}{2(\alpha-1)} ~\!y~\!\Delta_x\P^y_{\alpha-1}$ 
provided that $\alpha>1\,;$
\item[$iii)$]
$
\bDelta_b^m \P^y_\alpha=
\frac{\Gamma\big(\frac{1+b}{2}+\alpha\big)\Gamma(\alpha-m)}
{\Gamma\big(\frac{1+b}{2}+\alpha-m\big)\Gamma(\alpha)}~\!\Delta_x^m\P^y_{\alpha-m}
$~
for any positive integer $m<\alpha$.
\end{itemize}
\end{Lemma}

\proof
The  identities in $i)$  follow by direct computation, use the    identity 
$(n+2\alpha)c_{n,\alpha}=2\alpha c_{n,\alpha+1}$. 

If $\alpha>1$ we calculate
$$
\Delta_x |z|^{-n-2\alpha+2}=(n+2\alpha-2)\big(2\alpha|z|^{-n-2\alpha}-(n+2\alpha)y^2|z|^{-n-2\alpha-2}\big)\,.
$$
Since $\Delta_x\P^y_{\alpha-1}=c_{n, \alpha-1}y^{2\alpha-2}\Delta_x|z|^{-n-2\alpha+2}$ and $(n+2\alpha-2)c_{n,\alpha-1}=2(\alpha-1)c_{n,\alpha}$, we readily get
\begin{equation}
\begin{aligned}
\Delta_x\P^y_{\alpha-1}=&~2(\alpha-1)~\!y^{-2}
\big(2\alpha c_{n,\alpha} y^{2\alpha}|z|^{-n-2\alpha}-2\alpha c_{n,\alpha+1}y^{2\alpha+2}|z|^{-n-2\alpha-2}\big)\\
=&~
4\alpha(\alpha-1) ~\!y^{-2} \big(\P^y_\alpha-\P^y_{\alpha+1}\big)\,,
\label{eq:R1}
\end{aligned}
\end{equation}
which concludes the proof of $ii)$. 
If $m=1$ then the identity in $iii)$ follows from $i)$ and (\ref{eq:R1}). To conclude the proof of
$iii)$ use  induction.
\QED 

The next result deals with the  Fourier transform 
of the function $x\mapsto \P^y_\alpha(x)$ for $y\neq 0$.
In fact, $\widehat{\P^y_\alpha}(\xi)$ can be expressed via  Bessel functions, see for instance the computations in \cite{FrLe}. We provide a simple proof based on Lemma \ref{L:R1}.

\begin{Lemma}
\label{L:R2}
Let $\alpha>0$ and $y\neq 0$. Then 
\[
\widehat{\P^y_\alpha}(\xi)= 
\frac{2^{1-\alpha}}{(2\pi)^{\frac{n}{2}}\Gamma(\alpha)}~\!
|y|\xi||^\alpha K_\alpha(|y|\xi||)\,,
\]
where $K_\alpha$ is the (standard) modified Bessel function of the second kind of order $\alpha$.
\end{Lemma}

\proof
We can assume that $y>0$.
Since  $\P^y_\alpha$ is a radial function on $\R^n$,
$\|\P^y_\alpha\|_{L^1(\R^n)}= 1$ and 
 $\P^y_\alpha (x) = y^{-n} \P^1_\alpha\big(\frac{x}{y}\big)$, then 
$\widehat{\P^y_\alpha}(\xi)=\widehat{\P^y_\alpha}(|\xi|)$ is radial as well, 
$ (2\pi)^{n}\|\widehat{\P^y_\alpha}\|^2_\infty\le 1$ and
\begin{equation}
\label{eq:R3}
\widehat{\P^y_\alpha}(\rho)= \widehat{\P^1_\alpha}(\rho y)~\!,\qquad \rho=|\xi|~\!.
\end{equation}
Next, notice that  $\Delta \P^y_\alpha=(2\alpha-1)y^{-1}\partial_y \P^y_\alpha$ by $i)$ in Lemma \ref{L:R1}.
Hence, for $\rho>0$ fixed we have that the function $y\mapsto \widehat{\P^y_\alpha}=\widehat{\P^y_\alpha}(\rho)$ solves
\begin{equation}
\label{eq:R2}
(\widehat{\P^y_\alpha})''-(2\alpha-1)y^{-1}(\widehat{\P^y_\alpha})'-\rho^2 \widehat{\P^y_\alpha}=0
\end{equation}
on $\{y>0\}$.
We define
$K(\rho):=\rho^{-\alpha} \widehat{\P^1_\alpha}(\rho)$, 
 so that 
 $$
 \widehat{\P^y_\alpha}(\rho)=(\rho y)^\alpha  K(\rho y)
 $$ by (\ref{eq:R3}). Comparing with 
 (\ref{eq:R2}), we see that $K=K(t)$ solves
 $$
 t^2 K''+t K'-(t^2+\alpha^2)K=0~\!.
$$
It follows that $K$ is proportional to $K_\alpha$, the standard decreasing modified Bessel function of second kind of order $\alpha$. Since $\widehat{\P^y_\alpha}(0)=(2\pi)^{-n/2}$ and 
$t^\alpha K_\alpha(t)=2^{\alpha-1}\Gamma(\alpha)+o(1)$ as $t\to 0$, the proportionality constant is determined.
\QED
 
\subsection{Extension operators $\E_\alpha[u]=u*\P^y_\alpha$}
\label{SS:mia}

The next Lemma immediately implies Theorem \ref{T:mia}. 

\begin{Lemma}
 \label{L:L2_continuita}
 Let  $\alpha>0$, and assume that the pair $a,b\in \R$ satisfies
 $-1<b<1$, $0\le  2a < {n+1+b}$. There exists a constant $\gamma = \gamma(n,a,b)$ not depending on $\alpha$, such that 
$$
\iint\limits_{\R^{n+1}}|y|^{b}|z|^{-2a}|\E_\alpha[u]|^2~\!dz \le  \gamma \irn|x|^{-2a+b+1}|u|^2~\!dx
\quad\text{for any $u\in L^2(\R^n;|x|^{-2s}dx)$.}
$$
\end{Lemma}

\proof
By Remark \ref{R:M} we have that $|y|^b|z|^{-2a} \in A_2(\R^{n+1})$. 
Fix $u\in L^2(\R^n;|x|^{-2s}dx)$. Since $u\in L^1_{\rm loc}(\R^{n})$, then $\E_\alpha[u]=u*\P^y_\alpha$ is well defined 
and measurable on $\R^n$ for any $y\in \R$. 

Recall that $\|\P^y_\alpha\|_{1}=1$. Using the radial symmetry of the function $x\mapsto \P^y_\alpha(x)$, one can easily 
generalize \cite[Theorem III.2.2.(a)]{S}
to estimate
$$
|\E_{\alpha}[u](x, y)|=|u*\P^y_\alpha (x)| \leq 
|M_n[u](x)|\,, \quad \text{for almost every $(x,y) \in \R^{n+1}$}~\!,
$$
where $M_n$ is the operator in \eqref{eq:HLop}. We  identify $u=u(x)$ with a function $u=u(x,y)\in L^1_{\rm loc}(\R^{n+1})$
which is constant in the last variable. Thus
we can write
\[
|M_n[u](x)| = \sup_{r>0}\frac{1}{|B_r(x)|}\int\limits_{B_r(x)}|u| d\xi = 
\sup_{r>0} \frac{1}{2r|B_r(x)|}\int\limits_{y-r}^{y+r}d\tau\!\!\!\int\limits_{B_r(x)}|u|~\!d\xi \leq c(n)~\! |M_{n+1}[u](x,y)|\,,
\]
for any $(x,y)\in\R^{n+1}$, where $c(n)$ depends only on the dimension.
We infer that
\[
\iint\limits_{\R^{n+1}}|y|^b|z|^{-2a}|\E_{\alpha}[u]|^2~\!dz \leq c(n)~\!\iint\limits_{\R^{n+1}}|y|^b|z|^{-2a}|M_{n+1}[u]|^2~\!dz\le 
 c(n,a,b)~\!\iint\limits_{\R^{n+1}}|y|^b|z|^{-2a}|u(x)|^2~\!dz
\]
by the characterization of the Muckenhoupt class $A_2(\R^{n+1})$. The proof is complete.
\QED

\medskip
Next, we fix $u\in C^\infty_c(\R^n)$ and study
the regularity  of $\E_\alpha[u]$. Firstly, notice that $\E_\alpha[u]\in C^{0,\sigma}_{\me}(\R^{n+1})$ for any $\sigma\in(0,\min\{2\alpha,1\})$ and in particular
\begin{equation}
\label{eq:extreg}
 \lim\limits_{y\to 0}\|\E_{\alpha}[u](\cdot, y) - u\|_{\infty} = 0~\!.
\end{equation}
We recall the argument. Since $\|\P^y_\alpha\|_{L^1(\R^n)} = 1$, we can write
\[
|\E_{\alpha}[u] - u|
\leq c_{n,\alpha}\!\int\limits_{\{|x-\xi|\leq \delta\}}\! \frac{\|\nabla u\|_{\infty} |x-\xi|y^{2\alpha} }{(|x-\xi|^2 + |y|^2)^{\frac{n+2\alpha}{2}}}d\xi +c_{n,\alpha}\!\int\limits_{\{|x-\xi|>\delta\}} \!\frac{2\|u\|_{\infty}y^{2\alpha} }{(|x-\xi|^2 + |y|^2)^{\frac{n+2\alpha}{2}}} d\xi\,.
\]
The limit in (\ref{eq:extreg}) and the H\"older continuity of $\E_{\alpha}[u]$ follow since
$$
\begin{cases}
\|\E_{\alpha}[u](\cdot,y) - u\|_{\infty}=O(|y|+|y|^{2\alpha})&\text{if $\alpha\neq \frac12$}\\
\|\E_{\alpha}[u](\cdot,y) - u\|_{\infty}=O(|y|\log|y|)&\text{if $\alpha= \frac12$}
\end{cases}\qquad\text{as $y\to 0$.}
$$

In fact, the regularity of $\E_\alpha[u]$ increases for larger values of the parameter $\alpha$.

\begin{Lemma}
\label{L:extreg1}
Let $\alpha>1$ be not  integer, $u\in C^\infty_c(\R^n)$. 
Then 
\begin{eqnarray}
\label{eq:derivata1}
\partial^{2{m}}_y \E_{\alpha}[u] &=&\frac{\Gamma\big(\alpha+\frac12\big)}{\Gamma(\alpha)}
\sum_{\ell=0}^{m} \binom{{m}}{\ell}  
\frac{(-1)^\ell\Gamma(\alpha-\ell)}{\Gamma\big(\alpha+\frac12-\ell\big)} ~\E_{\alpha-\ell}[(-\Delta)^{m}u]
\\
\label{eq:derivata2}
y^{-1}\partial^{2{m}-1}_y \E_{\alpha}[u] &=&\frac{\Gamma\big(\alpha+\frac12\big)}{2\Gamma(\alpha)}
\sum_{\ell=0}^{{m}-1} \binom{{m}-1}{\ell}  
\frac{(-1)^{\ell+1}\Gamma(\alpha-\ell-1)}{\Gamma\big(\alpha+\frac12-\ell\big)} ~\E_{\alpha-\ell-1}[(-\Delta)^{{m}}u]
\end{eqnarray}
for  ${m}=1,\dots,[\alpha]$ and $y\neq 0$. Therefore, $\E_\alpha[u]\in C^{2[\alpha],\sigma}_{\me}(\R^{n+1})$ for any $\sigma\in(0,\min\{2\alpha,1\})$ and  the  Taylor expansion formula 
$$
\E_\alpha[u](\cdot,y)=\sum_{m=0}^{[\alpha]} \frac{\kappa_{\alpha,m}}{(2m)!}~\!y^{2m}~\!(-\Delta)^m u+o(y^{2[\alpha]})\quad \text{as $y\to 0$,}
$$
holds uniformly on $\R^n$, where
$$
\kappa_{\alpha,m}=
\frac{\Gamma\big(\alpha+\frac12\big)}{\Gamma(\alpha)}
\sum_{\ell=0}^m \binom{m}{\ell}  
\frac{(-1)^\ell\Gamma(\alpha-\ell)}{\Gamma\big(\alpha+\frac12-\ell\big)} \,.
$$
\end{Lemma}
\proof
By $i)$ (with $b=0$) and $ii)$ in Lemma \ref{L:R1} we have
$$
\Delta\P^y_\alpha=2\alpha(2\alpha-1)y^{-2}(\P^y_\alpha-\P^y_{\alpha+1})=
\frac{2\alpha-1}{2(\alpha-1)}\Delta_x\P^y_{\alpha-1}\,.
$$
Thus
\begin{equation}
\label{eq:ext2der}
\partial^2_y\E_\alpha[u]=(\Delta \E_\alpha[u]-\Delta_x\E_\alpha[u])=-
\frac{2\alpha-1}{2(\alpha-1)} \E_{\alpha-1}[-\Delta u] + \E_{\alpha}[-\Delta u]
\end{equation}
for $y\neq 0$. 
One can use induction and \eqref{eq:ext2der} to obtain (\ref{eq:derivata1}). Then (\ref{eq:derivata2}) follows,
since 
$$
\partial_y\E_{\alpha-\ell}[u]=u*\partial_y\P^y_{\alpha-\ell}=\frac{1}{2(\alpha-\ell-1)} y~\!(u*\Delta_x\P^y_{\alpha-\ell-1})=
-\frac{1}{2(\alpha-\ell-1)} y~\!\E_{\alpha-\ell-1}[-\Delta u]
$$
by $ii)$ in Lemma \ref{L:R1}.

We already observed that $\E_\alpha[u]$ is smooth
outside $\{y=0\}$ and that $\E_{\alpha-\ell}[(-\Delta)^m u]\in C^{0,\sigma}_{\me}(\R^{n+1})$ for any $m\ge 0$,
$\ell=0,\dots,[\alpha]$. Thus  
$\E_\alpha[u]\in C^{2[\alpha],\sigma}(\R^{n+1})$. The coefficients of Taylor formula for the even function $\E_\alpha[u]$
can be computed thanks to (\ref{eq:derivata1}), (\ref{eq:derivata2}), and taking (\ref{eq:extreg}) into account
 (with $\alpha, u$ replaced by 
$\alpha-\ell, (-\Delta)^mu$, respectively).
\QED

\medskip

We conclude this section by studying the behaviour of $\E_\alpha$ acting on smooth functions and then
on the space 
$\mathcal D^{s}(\R^n)$. Recall that $\E_{\alpha}[u] \in C^{2[\alpha],\sigma}_\me(\R^{n+1})$ for $u\in C^\infty_c(\R^n)$
by Lemma \ref{L:R2}. 

\begin{Lemma}\label{L:Walphasumm}
Let $\alpha>[\alpha]\ge [s]\ge 1$, $b=1-2(s-[s])$. There exists a constant 
$C_\alpha$ depending only on $n, s$ and $\alpha$, such that
$$
\iint\limits_{\R^{n+1}}|y|^b|\bnabla_{\!b}^{1+[s]}\E_\alpha[u]^2~\!dz=  C_\alpha \irn|\Dshalf u|^2~\!dx
\quad\text{ for any $u\in C^\infty_c(\R^n)$. }
$$
\end{Lemma}

\proof 
Fix $u\in C^\infty_c(\R^n)$. If $[s]= 2m-1$ is odd, we have

\[
 \bnabla^{1+[s]}_b \E_{\alpha}[u] = \bDelta^{m}_b \E_{\alpha}[u] = u*(\bDelta^{m}_b \P^y_{\alpha}) = c \, u*\Delta_x^{m}\P^y_{\alpha-m} = c~(\Delta^{\!m} u)*\P^y_{\alpha-m}~\!,
\]
by $iii)$ in Lemma \ref{L:R1} and since $\alpha>[s]\ge m$.
Thus, using also  Lemma \ref{L:R2} we infer
 \begin{eqnarray}
 \nonumber
\iint\limits_{\R^{n+1}} |y|^b| \bnabla^{1+[s]}_b \E_{\alpha}[u]|^2~\!dz
&=&c \irn  |\xi|^{2(m+\alpha)} |\hat{u}|^2~\!d\xi 
\int\limits_{0}^\infty y^{b+2(\alpha-m)} |K_{\alpha-m}(y|\xi|)|^2~\!dy\\
&=&
c \irn  |\xi|^{4m-b-1} |\hat{u}|^2~\!d\xi 
\int\limits_{0}^\infty t^{b+2(\alpha-m)} |K_{\alpha-m}|^2~\!dt\,.
 \label{eq:R5}
 \end{eqnarray}
Recall that $b=1-2(s-[s])=-2s +(4m-1)$, so that 
$\irn  |\xi|^{4m-b-1} |\hat{u}|^2~\!d\xi =\irn |\Dshalf u|^2~\!dx$. Since $b>-1$, $t^{\alpha-m}K_{\alpha-m}(t)=O(1)$ as $t\to 0$ and $K_{\alpha-m}(t)$ decays exponentially
as $t\to \infty$, then the last integral in (\ref{eq:R5}) converges, and the "odd" case is over.

\medskip
If $[s]=2m\ge 2$ is even, then $\alpha > m+1$. With similar computations we find
\begin{multline*}
\iint\limits_{\R^{n+1}}|y|^b |\bnabla^{1+[s]}_b \E_{\alpha}[u]|^2\,~\!dz = \iint\limits_{\R^{n+1}}|y|^b |\nabla_{\!x}\bDelta^m_b \E_{\alpha}[u]|^2\,~\!dz\, + \iint\limits_{\R^{n+1}}|y|^b |\partial_y\bDelta^m_b \E_{\alpha}[u]|^2\,~\!dz\,\\
 = c\iint\limits_{\R^{n+1}}|y|^b|\nabla_x \big((\Delta^m u) * \P^y_{\alpha-m}\big)|^2~\!dz+ c\iint\limits_{\R^{n+1}}|y|^{b+2}|(\Delta^{m+1} u) * \P^y_{\alpha-m-1}|^2~\!dz\\
=c\int\limits_{\R^n}|\xi|^{4m+1-b}|\hat u|^2d\xi\int\limits_{0}^\infty t^{b+2(\alpha-m)}(|K_{\alpha-m}|^2 +|K_{\alpha-m-1}|^2)  dt\,. 
\end{multline*}
The proof  is concluded. \QED

\begin{Lemma}
\label{L:semplice}
Let $s>1$ be not an integer, $\alpha>[\alpha]\ge [s]\geq 1$, $b=1-2(s-[s])$. Then
\begin{itemize}
\item[$i)$] $\E_\alpha[u]\in \mathcal D^{1+[s];b}_\me(\R^{n+1})$ for any $u\in \mathcal D^s(\R^n)$;
\item[$ii)$] $\E_\alpha:\mathcal D^s(\R^n)\to \mathcal D^{1+[s];b}_\me(\R^{n+1})$ is, up to a constant, an isometry;
\item[$iii)$] $\textrm{\rm Tr}(\E_\alpha[u])=u$ for any $u\in \mathcal D^s(\R^n)$.
\end{itemize}
\end{Lemma}

\proof
The proof $i)$ is quite technical and it is postponed to the Appendix. Claim $ii)$ is an immediate consequence of
Lemma \ref{L:Walphasumm}. Since $\text{Tr}(\E_\alpha[u])=u$ if $u\in C^\infty_{c}(\R^n)$, then $iii)$ follows from
$ii)$, as $C^\infty_{c}(\R^n)$ is dense in $\mathcal D^s(\R^n)$.
\QED

\section{Proof of Theorem \ref{T:main}}
\label{S:proof}

If $s\in (0,1)$, then Theorem \ref{T:main} follows by adapting the proofs in \cite{FrLe} (see also 
\cite[Section 5]{FLS}). Therefore, from now on we assume that $s>1$.

We start by writing $iii)$ in Lemma \ref{L:R1} with $\alpha=s$ and $b=1-2(s-[s])$. We have

\begin{equation}
\label{eq:R8}
(-\Delta_x)^{\nu}\P^y_{s-{\nu}}= 
\frac{([s]-{\nu})!}{[s]!}~\frac{\Gamma(s)}{\Gamma(s-{\nu})}~\!(-\bDelta_b)^{\nu} \P^y_s~,
\qquad {\nu}=1,\dots,[s].
\end{equation}

Clearly $i)$ follows from Lemma \ref{L:semplice} with $\alpha=s$. 

Thanks to the continuity of the map $\E_s:\mathcal D^s(\R^n)\to\mathcal D^{1+[s];b}_\me(\R^{n+1})$,
 in order to prove $ii)$ and $iii)$ we can assume that $u\in C^\infty_c(\R^n)$. Then Lemma \ref{L:extreg1} gives
$$
\E_s[u]\in C^{2[s],\sigma}(\R^{n+1})\quad \text{for any $\sigma\in (0,1)$}~\!.
$$
We fix  ${V}\in C^\infty_{c;\me}(\R^{n+1})$ and use (\ref{eq:R8}) with $\nu=[s]$ to compute 
$$
((-\Delta)^{[s]}u)*\P^y_{s-[s]}=u*((-\Delta_x)^{[s]}\P^y_{s-[s]})=\frac{\Gamma(s)}{[s]!\Gamma(s-[s])}
~\!u*((-\bDelta_b)^{[s]}\P^y_{s})=\frac{d_{s-[s]}}{d_{s}}(-\bDelta_b)^{[s]}(u*\P^y_s)~\!.
$$
We see that
$$
\E_{s-[s]}[(-\Delta)^{[s]} u]=\frac{d_{s-[s]}}{d_s}(-\bDelta_b)^{[s]}\E_s[u]\in \mathcal D^{1;b}_{\me}(\R^{n+1})
$$
because $s-[s]\in(0,1)$ and thanks to the results in \cite{CS, FrLe, FLS}. Moreover,
$$
\begin{aligned}
2d_{s-[s]}\langle \Ds u,\textrm{Tr}({V})\rangle&=
2d_{s-[s]}\langle (-\Delta)^{s-[s]}(-\Delta)^{[s]} u,\textrm{Tr}({V})\rangle\\
&=
\iint\limits_{\R^{n+1}} |y|^b\nabla  \E_{s-[s]}[(-\Delta)^{[s]}u]\cdot \nabla V=
\frac{d_{s-[s]}}{d_s} \iint\limits_{\R^{n+1}} |y|^b\nabla(-\bDelta_b)^{[s]}\E_s[u]\cdot \nabla V.
\end{aligned}
$$
Thus integration by parts and Lemma \ref{L:Step1} (with $W=\E_s[u]$ and $k=1+[s]$) give
$$
2d_s\langle \Ds u,\textrm{Tr}({V})\rangle =
\iint\limits_{\R^{n+1}} |y|^b(-\bDelta_b)^{[s]}\E_s[u] (-\bDelta_b) V=
\iint\limits_{\R^{n+1}} |y|^b\bnabla_{\!b}^{1+[s]}\E_s[u]~\!\bnabla_{\!b}^{1+[s]} V~\!.
$$
Since $C^\infty_{c;\me}(\R^{n+1})$ is dense in $\mathcal D^{1+[s];b}_\me(\R^{n+1})$, we have proved that 
$\E_s[u]$ satisfies (\ref{eq:aggiunto}).
In fact, the variational problem
$$
\begin{cases}
~U\in \mathcal D^{1+[s];b}_\me(\R^{n+1})\\
\displaystyle{~\iirn|y|^b\bnabla_{\!b}^{1+[s]} U~\bnabla_{\!b}^{1+[s]}V~\!dz=2d_s\langle\Ds u,\textrm{\rm Tr}(V)\rangle~~~
\text{for any $V\in \mathcal D^{1+[s];b}_\me(\R^{n+1})$}}
\end{cases}
$$
is equivalent to 
the minimization problem in  (\ref{eq:min1}), which evidently admits a unique solution. Thanks to (\ref{eq:aggiunto}),
we infer that
$\E_s[u]$ achieves the minimum in (\ref{eq:min1}), which concludes the proof of $ii)$.

\medskip

Next, the convex minimization problem  (\ref{eq:min2}) has a unique solution $U_0\in \mathcal D^{1+[s];b}_\me(\R^{n+1})$, which is 
the minimal distance projection of
$0\in \mathcal D^{1+[s];b}_\me(\R^{n+1})$ on the closed, affine space $\text{Tr}^{-1}\{u\}$. Thus, $U_0$ is the unique
point in $\text{Tr}^{-1}\{u\}$ which is orthogonal to $\text{Tr}^{-1}\{0\}$. Since (\ref{eq:aggiunto}) 
implies
$$
\iirn|y|^b\bnabla_{\!b}^{1+[s]} \E_s[u]~\bnabla_{\!b}^{1+[s]}V=0\quad 
\text{for any $V\in \text{Tr}^{-1}\{0\}$,}
$$
we can conclude that  $U_0=\E_s[u]$.

\medskip

The proof of (\ref{eq:uguale_operatori}) is an adaptation of \cite[Proposition 3.5]{FrLe}, where $n=1$ and $s \in (0,1)$ are assumed. 
It suffices to study the behaviour of $y^b\partial_y (-\bDelta_b)^{[s]}\E_s[u]$ as $y\to 0^+$. For any fixed $y>0$
we compute the Fourier transform of 
\[
\mathbb L^y_s[u] :=  (-\Delta)^s u+\frac{1}{d_s}|y|^{b-1} y \partial_y (-\bDelta_b)^{[s]}\E_s[u]~\!.
\]
By (\ref{eq:R8})  we have 
$\partial_y(-\bDelta_b)^{[s]}\E_s[u](\cdot,y) = u *( \partial_y(-\bDelta_b)^{[s]} \P^y_s) = \frac{d_s}{d_{s-[s]}} u * ((-\Delta_x)^{[s]}\partial_y\P^y_{s-[s]})$. 
Thus
\[
\widehat {\mathbb L^y_s[u]} =  |\xi|^{2s}\hat u\ \Big(1+\frac{(2\pi)^{\frac{n}{2}}}{d_{s-[s]}}~\!y~\!(y|\xi|)^{b-1}\widehat{\partial_y \P_{s-[s]}^y}\Big)~\!.
\]
Put $\alpha:=s-[s]$, recall (\ref{eq:ds}) and notice that $b-1=-2\alpha$. Thanks to Lemma \ref{L:R2}, we infer that
$$
\widehat {\mathbb L^y_s[u]} = |\xi|^{2s}\hat u\Big(1+\frac{2^\alpha}{\Gamma(1-\alpha)}~~\!(y|\xi|)^{1-2\alpha}
\partial_t\big(t^\alpha K_\alpha(t)\big)_{|t=y|\xi|}\Big)
=: |\xi|^{2s}\hat u~\! \Phi_\alpha(y|\xi|).
$$
Using the known formulae for the modified Bessel functions, one can compute
$$
t^{1-2\alpha}\partial_t(t^\alpha K_\alpha(t))=t^{1-\alpha}(K'_\alpha(t)+\alpha t^{-1}K_\alpha(t))=
- t^{1-\alpha}K_{1-\alpha}(t),
$$
so that
$$
\Phi_\alpha(t)= 1-\frac{2^\alpha}{\Gamma(1-\alpha)}~\!t^{1-\alpha}K_{1-\alpha}(t)=o(1)\quad\text{as $t\to 0$}.
$$
It follows that $\Phi_\alpha(y|\xi|)\to 0$ almost everywhere and in the weak$^*$ topology of $L^\infty(\R^n)$ as $y\to 0$
(recall that $K_{1-\alpha}$ decays exponentially at infinity). 
Thus, for any $v\in \mathcal D^s(\R^n)$ we have
$$
\begin{aligned}
|\langle \mathbb L^y_s[u], v\rangle|^2 = \Big| \int\limits_{\R^n} |\xi|^{2s}\hat u~\!\overline{\hat v} \Phi_\alpha(y|\xi|)~\!d\xi\Big|^2
&\le  \int\limits_{\R^n} |\xi|^{2s}|{\hat v}|^2~\!d\xi\, \int\limits_{\R^n} |\xi|^{2s}|\hat u|^2 \Phi_\alpha(y|\xi|)^2~\!d\xi \\
 &= \|\Dshalf v\|_2^2  \int\limits_{\R^n} |\xi|^{2s}|\hat u|^2 \Phi_\alpha(y|\xi|)^2~\!d\xi
 \end{aligned}
$$
which, together with $|\xi|^{2s}|\hat u|^2\in L^1(\R^n)$, implies
$$
\|\mathbb L_s[u]\|^2_{\mathcal D^{-s}(\R^n)}\le  \int\limits_{\R^n} |\xi|^{2s}|\hat u|^2 \Phi_\alpha(y|\xi|)^2~\!d\xi=o(1)~\!.
$$
The proof of (\ref{eq:uguale_operatori}) is complete.

\medskip 
To prove (\ref{eq:andata}) we  use $ii)$ in Lemma \ref{L:R1} to obtain
$\Delta_x\P^y_{s-m}=2(s-m)y^{-1}\partial_y\P^y_{s-m+1}$. Thus
$(-\Delta_x)^m\P^y_{s-m}=-2(s-m)y^{-1}\partial_y(-\Delta_x)^{m-1}\P^y_{s-(m-1)}$.
By applying (\ref{eq:R8}) two times (with $\nu=m-1$ and then
$\nu=m$) we easily infer 
$(-\bDelta_b)^{m}\P^y_s = -2(1+[s]-m)~y^{-1}\partial_y(-\bDelta_b)^{m-1}\P^y_s$, 
and  (\ref{eq:andata})   follows.

\medskip
Now we prove (\ref{eq:limits}). If $s>2m$ then $(-\Delta)^mu\in\mathcal D^{s-2m}(\R^n)$. Using (\ref{eq:R8}) as before, we obtain 
$$
(-\bDelta_b)^m\E_s[u]= \frac{d_{s}}{d_{s-m}}\E_{s-m}[(-\Delta)^mu]\in {\mathcal D^{1+[s]-2m;b}_\me(\R^{n+1})}
$$
by Lemma \ref{L:semplice}. Theorem \ref{T:traces} applies and gives (\ref{eq:limits}), because
$$
\textrm{Tr}((-\bDelta_b)^m\E_s[u])= \frac{d_{s}}{d_{s-m}} \textrm{Tr}(\E_{s-m}[(-\Delta)^mu])= \frac{d_{s}}{d_{s-m}}(-\Delta)^mu.
$$

If $s<2m$ the proof can be obtained by repeating the argument for (\ref{eq:uguale_operatori}).
For $y>0$ we  formally define the operator $\mathbb L^y_s[u]$ in the dual space $\mathcal D^{s-2m}(\R^n)=\mathcal D^{2m-s}(\R^n)'$   by
$$
\mathbb L^y_s[u] := (-\Delta)^mu-\frac{d_{s-m}}{d_{s}}(-\bDelta)^m\E_s[u](\cdot,y)=
(-\Delta)^mu-\E_{s-m}[(-\Delta)^mu](\cdot,y).
$$
We have
$\widehat {\mathbb L^y_s[u]} =|\xi|^{2m}\hat{u}\big(1-(2\pi)^{n/2} \hat{\P}^y_{s-m}(\xi)\big)=|\xi|^{2m}\hat{u}~\!\Phi_\alpha(y|\xi|)$, where now $\alpha = s-m$ and
$$
\Phi_\alpha(y|\xi|)=1-\frac{2^{1-\alpha}}{\Gamma(\alpha)}~\!(y|\xi|)^\alpha K_\alpha(y|\xi|)\to 0\quad
\text{weakly$^*$ in $L^\infty(\R^n)$, as $y\to0^+$.}
$$
The conclusion follows as for (\ref{eq:uguale_operatori}). The limits in (\ref{eq:limits2}) can be checked in a similar way
via \eqref{eq:derivata1}, \eqref{eq:derivata2}, using Lemma \ref{L:semplice} and Theorem \ref{T:traces} if $s>2m$, and thanks to 
Lemma \ref{L:R1}, Lemma \ref{L:R2}, known properties of Bessel's functions if $s<2m$.

\medskip

The last assertion in Theorem \ref{T:main} readily follows from Lemma  \ref{L:extreg1}.
\QED
 
\subsection{On the variational problem (\ref{eq:0_extequation})}
\label{SS:distribuzionale}

Here we assume that $s\in(0,n/2)$ is not an integer and put $b=1-2(s-[s])$, as in Theorem \ref{T:main}.

Let $U\in C^\infty_{c;\me}(\R^{n+1})$. Then 
$-\div(|y|^b\nabla (-\bDelta_b)^{[s]} U)=|y|^b(-\bDelta_b)^{1+[s]}U\in L^1_{\rm loc}(\R^{n+1})$ can be regarded
as a distribution on $\R^{n+1}$ which vanishes on functions that are odd in the $y$-variable. 

In addition, for any $\f\in C^\infty_{c;\me}(\R^{n+1})$ we can integrate by parts to get
$$
\begin{aligned}
 \langle -\div(|y|^b\nabla  (-\bDelta_b)^{[s]} U),\f\rangle&=\iirn |y|^b\nabla  (-\bDelta_b)^{[s]} U\cdot\nabla \f~\!dz=
 \iirn |y|^b  (-\bDelta_b)^{[s]} U\cdot(-\bDelta_b)\f~\!dz\\
& = \iint\limits_{\R^{n+1}} |y|^b\bnabla_{\!b}^{1+[s]}U\, \bnabla_{\!b}^{1+[s]} \f~\!dz
 \end{aligned}
 $$
 by Lemma \ref{L:Step1}. 
Since $C^\infty_{c,\me}(\R^{n+1})$ is dense in $\mathcal D^{1+[s];b}_\me(\R^{n+1})$, we see that for any 
$U\in \mathcal D^{1+[s];b}_\me(\R^{n+1})$ we can look at
$-\div(|y|^b\nabla  (-\bDelta_b)^{[s]} U)$ as the distribution in the dual space $\mathcal D^{1+[s];b}_\me(\R^{n+1})'$
which acts as follows,
$$
\langle -\div(|y|^b\nabla  (-\bDelta_b)^{[s]} U), V\rangle= \iint\limits_{\R^{n+1}} |y|^b\bnabla_{\!b}^{1+[s]}U \bnabla_{\!b}^{1+[s]} V~\!dz
\quad \text{for any $V\in \mathcal D^{1+[s];b}_\me(\R^{n+1})$}~\!.
$$

\medskip

Next,  let $u\in \mathcal D^s(\R^n)$. The trace map $\text{Tr}: {\mathcal D}_\me^{1+[s];b}(\R^{n+1})\to \mathcal D^s(\R^n)$ in Theorem \ref{T:traces} can be composed
with $\Ds u$, which is a linear form on $\mathcal D^s(\R^n)$. Instead of writing $(\Ds u)\circ \text{Tr}$, we 
prefer to use the more suggestive notation $\delta_{\{y=0\}}\Ds u$. Thus
$$
\delta_{\{y=0\}}\Ds u \in {\mathcal D}_\me^{1+[s];b}(\R^{n+1})',~~ 
\langle \delta_{\{y=0\}}\Ds u,V \rangle= \langle \Ds u, \text{Tr}(V)\rangle\quad \text{for any $V\in \mathcal D^{1+[s];b}_\me(\R^{n+1})$.}
$$

In conclusion, we gave a precise interpretation of the differential equation in (\ref{eq:0_extequation}) as an
equality in the dual space ${\mathcal D}_\me^{1+[s];b}(\R^{n+1})'$. Moreover, (\ref{eq:0_extequation}) gives
the Euler-Lagrange equations for the minimization problem (\ref{eq:min1}).

\appendix

\section{\!\!\!\!\!\!ppendix: proof of $i)$ in Lemma \ref{L:semplice}. 
}

Thanks to Lemma \ref{L:Walphasumm}, we only need to show that 
$\E_\alpha[u]\in \D^{1+[s];b}_{\me}(\R^{n+1})$ for a fixed $u\in C^\infty_c(\R^n)$.

The idea is quite simple. We take a cut-off function $\f\in C^\infty_{c;\me}(\R^{n+1})$ such that $\f\equiv 1$ in a neighbourhood of 
the origin and put $\f_\lambda(z)=\f(\lambda^{-1}z)$, $\lambda >0$. Then 
$$
\f_\lambda \E_\alpha[u]\in C^{1+[s]}_{c;\me}(\R^{n+1}) \subset \D^{1+[s];b}_\me(\R^{n+1}),
$$
by  Lemma \ref{L:extreg1} and Remark \ref{R:ck_dense}. Evidently,
$\f_\lambda \E_\alpha[u]\to \E_\alpha[u]$ in $L^2(\R^{n+1}; |y|^b|z|^{-2(1+[s])})$ and almost everywhere as $\lambda \to \infty$. 
To conclude the proof we will show that
\begin{equation}
\label{eq:crucial}
\|\f_\lambda\E_\alpha[u]\|_{1+[s]; b} \leq c\|\E_\alpha[u]\|_{1+[s];b}\,,
\end{equation}
which implies that $\f_\lambda \E_\alpha[u]\to \E_\alpha[u]$ weakly in $\D^{1+[s];b}_{\me}(\R^{n+1})$\footnote{\footnotesize{in fact, $\f_\lambda \E_\alpha[u]\to \E_\alpha[u]$ in the $\D^{1+[s];b}_{\me}(\R^{n+1})$-norm.
In order to skip quite long computations, we limit ourselves to prove the weak convergence, which
is enough for our purposes.}}.

From now on we neglect to write the volume integration forms $dz$ on $\R^{n+1}$ and $dx$ on $\R^n$. Also,
accordingly with (\ref{eq:nn}) (for $a=0$), we put
$$
\nn U\nn_{k;b}^2=\nn U\nn_{k;0,b}^2 := 
\sum_{j=0}^{k}\ \iint\limits_{\R^{n+1}}|y|^b|z|^{-2(k-j)}|\bnabla_{\!b}^j U|^2\,.
$$
We divide the remaining part of the proof in two steps.

\medskip

\paragraph{Step 1: an estimate.}
First of all we prove the inequality
\begin{equation}
\label{eq:all}
\nn \E_{\alpha}[v]\nn_{1+[s];b}\le c \|\E_\alpha[v]\|_{1+[s]; b}\,, \quad \text{ for any $v \in C^\infty_c(\R^n)$,} 
\end{equation}
(which, of course, can not be derived via Theorem \ref{T:Hardy}).
Thanks to Lemma \ref{L:Walphasumm}, it is 
sufficient to prove
\[
\iint\limits_{\R^{n+1}}|y|^{b}|z|^{-2(1+[s]-j)}|\bnabla_{\!b}^{j}\E_{\alpha}[v]|^2 \leq c \irn |\Dshalf  v|^2\,,\quad j=0,\dots  [s].
\]
We start by noticing that $iii)$ in Lemma \ref{L:R1} implies
\begin{equation}
\label{eq:Wbdeltam}
\bDelta_b^m \E_{\alpha}[v] = 
c~\! \E_{\alpha-m}[ (-\Delta)^mv] \qquad \text{for any integer $0\le m<\alpha$} ~\!,
\end{equation}
where the constant $c=c(s,\alpha,m)$ does not depend on $v$.

If $j = 2m$ is even, then $\bnabla_{\!b}^{j}\E_{\alpha}[v] = \bDelta_b^m \E_{\alpha}[v]$.
Thus, Lemma \ref{L:L2_continuita} with $s$, $ \alpha$ and $u$ replaced by $s-2m$, $\alpha-m$ and $(-\Delta)^m u$, respectively,
gives
$$
\iint\limits_{\R^{n+1}}\!\!|y|^{b}|z|^{-2(1+[s]-j)}|\bnabla_{\!b}^{j}\E_{\alpha}[v]|^2= c\!\!\iint\limits_{\R^{n+1}}\!\!|y|^{b}|z|^{-2(1+[s]-2m)}|\E_{\alpha-m}[ (-\Delta)^mv]|^2 \\
\le c\! \int\limits_{\R^n}\!\!|x|^{-2(s-2m)}|(-\Delta)^m v|^2\,.
$$
Next we use the Hardy inequality for the fractional Laplacian $(-\Delta)^{s-2m}$ to infer
$$
\iint\limits_{\R^{n+1}}|y|^{b}|z|^{-2(1+[s]-j)}|\bnabla_{\!b}^{j}\E_{\alpha}[v]|^2
\le c \irn |(-\Delta)^{\!\frac{s-2m}{2}} (-\Delta)^m v|^2
=  c \irn |\Dshalf  v|^2~\!.
$$

{If  $j= 2m+1$ is odd, we write
\[
\begin{aligned}
\big|\bnabla_{\!b}^{j}\E_{\alpha}[v]\big|^2 = c\sum_{\ell=1}^n
\big|\E_{\alpha-m}[(-\Delta)^m\partial_{x_\ell}v]\big|^2
+c\big|\partial_y\E_{\alpha-m}[(-\Delta)^m v]\big|^2
\,.
\end{aligned}
\] 
We use again Lemma \ref{L:L2_continuita} (for the exponents 
$s-2m-1, \alpha-m$) and then the Hardy inequality for the fractional Laplacian $(-\Delta)^{s-2m-1}$ to get
\begin{multline*}
\sum_{\ell=1}^n\ \iint\limits_{\R^{n+1}}|y|^{b}|z|^{-2(1+[s]-2m-1)}|\E_{\alpha-m}[(-\Delta)^m\partial_{x_\ell}v]|^2
\leq 
c\sum_{\ell=1}^n \ \int\limits_{\R^n}|x|^{-2(s-2m-1)}|(-\Delta)^m \partial_{x_\ell}v|^2 \\
\le 
c \sum_{\ell=1}^n\ \irn |(-\Delta)^{\frac{s-1}{2}}  \partial_{x_\ell}v|^2 =c \irn |(-\Delta)^{\frac{s}{2}}v|^2~\!.
\end{multline*}

To handle the weighted $L^2$ norm 
of $\partial_y\E_{\alpha-m}[(-\Delta)^m v]$ we first
use $i)$ in Lemma \ref{L:R1} to get 
\[
|\partial_y\E_{\alpha-m}[(-\Delta)^m v]|=
c|y\bDelta_b\E_{\alpha-m}[(-\Delta)^{m} v]|
\le c|z|~\!\big|\bDelta_b\E_{\alpha-m}[(-\Delta)^{m} v]\big|\,,
\]
and next the equality $\bDelta_b\E_{\alpha-m}[ (-\Delta)^m v]=c
\bDelta_b^{m+1}\E_{\alpha}[ v]$ 
(compare with \eqref{eq:Wbdeltam}). 
Since $2m+1=j$ and $\bDelta_b^{m+1}=\bnabla_{\!b}^{2m+2}=\bnabla_{\!b}^{1+j}$, we find 
\[
\iint\limits_{\R^{n+1}}|y|^{b}|z|^{-2(1+[s]-2m-1)}
|\partial_y\E_{\alpha-m}[(-\Delta)^m v]|^2 \leq c \iint\limits_{\R^{n+1}}|y|^{b}|z|^{-2(1+[s]-(1+j))}
|\bnabla_{\!b}^{1+j}\E_{\alpha}[v]|^2\,.
\]

Let $j < [s]$. Since $1+j$ is even and $1+j \leq [s]$, we can argue as in the ''even'' case. We obtain
\[
\iint\limits_{\R^{n+1}}|y|^{b}|z|^{-2(1+[s]-(1+j))}
|\bnabla_{\!b}^{1+j}\E_{\alpha}[v]|^2 \leq c\int\limits_{\R^n}|x|^{-2(s-(1+j))}|(-\Delta)^{\frac{1+j}{2}}v|^2   \leq c\irn |(-\Delta)^{\frac{s}{2}}v|^2\,.
\]
If $j= [s]$ we have that
\[
\iint\limits_{\R^{n+1}}|y|^{b}|z|^{-2(1+[s]-(1+j))}
|\bnabla_{\!b}^{1+j}\E_{\alpha}[v]|^2 =\iint\limits_{\R^{n+1}}|y|^{b}|\bnabla_{\!b}^{1+[s]}\E_{\alpha}[v]|^2
=c \irn|\Dshalf v|^2~\!
\]
by Lemma \ref{L:Walphasumm}.
The proof of (\ref{eq:all}) is complete.

\medskip
\paragraph{Step 2: proof of (\ref{eq:crucial}).}
~We introduce the sets
\[
\mathcal{A_\beta}:= \{\ \phi \in C^\infty_{c;\me}(\R^{n+1}) \ |\ \phi(z) =
\beta \text{ in a neighbourhood of $0$}
\ \}~,\quad \beta=0,1\,,
\]
so that $\f\in \mathcal A_1$, while the partial derivatives of $\f$  of any order belong to $\mathcal A_0$.
For any $\phi \in \mathcal{A}_\beta, \lambda>0$, we put $\phi_\lambda(z) = \phi(\lambda^{-1}z)$.
By direct computation one gets, for any $j, m \geq 0$ integers, 
\begin{gather}
\label{eq:lambdaA}
\partial_{x_\ell}\bDelta^m_b\phi_\lambda =\lambda^{-2m-1} (\partial_{x_\ell}\bDelta^m_b\phi)_\lambda, \quad\quad y \partial_y \bDelta^m_b\phi_\lambda = \lambda^{-2m}(y \partial_y \bDelta^m_b\phi)_\lambda\,;\\
\label{eq:lambdaB}
\displaystyle |\bnabla_{\!b}^j \phi_\lambda |  \leq c\|\bnabla_{\!b}^j \phi \|_\infty |z|^{-j}\,. 
\end{gather}

Next we prove, by induction on $k$, the crucial estimate 
\begin{equation}
\label{eq:crucial2}
\|\phi_\lambda\E_\alpha[v]\|_{1+k; b} \leq c\|\E_\alpha[v]\|_{1+k; b}\,,\qquad k:=[s]\ge 1~, \text{ for any $v \in C^\infty_c(\R^n)$, $\phi \in \mathcal{A}_\beta$.}
\end{equation}
Let $k=1$. We compute
\[
\bDelta_b(\phi_\lambda\E_\alpha[v]) = \phi_\lambda\bDelta_b \E_\alpha[v] + 2\nabla \E_\alpha[v]\cdot \nabla \phi_\lambda + \E_\alpha[v]\bDelta_b\phi_\lambda\,.
\]
Thus, using also  (\ref{eq:lambdaB}) and \eqref{eq:all} we obtain
\[
\|\phi_\lambda\E_\alpha[v]\|^2_{2; b} = \iint\limits_{\R^{n+1}}|y|^b|\bDelta_b(\phi_\lambda\E_\alpha[v])|^2
\leq c\sum_{j=0}^{2}\ \iint\limits_{\R^{n+1}}|z|^{-2(2-j)}|\bnabla^j_b \E_\alpha[v]|^2= c \nn \E_\alpha[v]\nn^2_{2; b} \leq c\| \E_\alpha[v]\|^2_{2; b} \,,
\]
where the constants $c >0$ do not depend on $\lambda$. Hence (\ref{eq:crucial2}) is proved in this case.

Next, fix $k \geq 2$ and assume that (\ref{eq:crucial2}) holds true for every $1 \leq k' \leq k-1$. 
Arguing by induction one can prove the  Leibniz-type formula
\[
\text{}\bnabla^{1+k}_b(\phi_\lambda\E_\alpha[v]) = \sum_{i=0}^{1+k}C_{i,k}\bnabla_{\!b}^i \E_\alpha[v] \bnabla_{\!b}^{1+k-i}\phi_\lambda + \sum_{j,h \text{\ odd}\atop{2 \leq j+h \leq k}}C_{j,h,k} \bnabla_{\!b}^{1+k-(j+h)}\big(\bnabla_{\!b}^j \E_\alpha[v] \cdot \bnabla_{\!b}^h \phi_\lambda \big)\,, 
\]
where $C_{i,k}, C_{j,h,k}>0$ only depend on $i$, $j$, $k$. Thus
\[
\|\phi_\lambda\E_\alpha[v]\|^2_{1+k; b} \leq c\sum_{i=0}^{1+k}\ \iint\limits_{\R^{n+1}}|y|^b|\bnabla_{\!b}^i \E_\alpha[v] \bnabla_{\!b}^{1+k-i}\phi_\lambda|^2+\, c\!\!\sum_{j,h \text{\ odd}\atop{2 \leq j+h \leq k}}\ \iint\limits_{\R^{n+1}}|y|^b |\bnabla_{\!b}^{1+k-(j+h)}\big(\bnabla_{\!b}^j \E_\alpha[v] \cdot \bnabla_{\!b}^h \phi_\lambda \big)|^2\,.
\]
Thanks to (\ref{eq:lambdaB}) and \eqref{eq:all} we readily get 
\[
\sum_{i=0}^{1+k}\ \iint\limits_{\R^{n+1}}|y|^b|\bnabla_{\!b}^i \E_\alpha[v] \bnabla_{\!b}^{1+k-i}\phi_\lambda|^2 \leq c \nn \E_\alpha[v]\nn_{1+k;b}^2 \leq c \| \E_\alpha[v]\|_{1+k;b}^2\,.
\]

To complete the inductive step and thus the proof of \eqref{eq:crucial2}, we show that for any couple of odd integers $j,h$ such that $2 \leq j+h \leq k$, it holds 
\begin{equation}\label{eq:tech03}
\iint\limits_{\R^{n+1}}|y|^b |\bnabla_{\!b}^{1+k-(j+h)}\big(\bnabla_{\!b}^j \E_\alpha[v] \cdot \bnabla_{\!b}^h \phi_\lambda \big)|^2 \leq c \|\E_\alpha[v] \|_{1+k;b}^2\,, 
\end{equation}
for some $c >0$ which does not depend on $\lambda$. We start the proof of \eqref{eq:tech03} by noticing
that
\[
\begin{aligned}
\bnabla_{\!b}^j \E_\alpha[v] \cdot \bnabla_{\!b}^h \phi_\lambda = \nabla\bDelta_b^{\frac{j-1}{2}} \E_\alpha[v] \cdot \nabla\bDelta_b^{\frac{h-1}{2}} \phi_\lambda &=  \sum_{\ell=1}^n\partial_{x_\ell}\bDelta_b^{\frac{j-1}{2}} \E_\alpha[v] \partial_{x_\ell}\bDelta_b^{\frac{h-1}{2}}\phi_\lambda\\
 &+\big(y^{-1}\partial_y\bDelta_b^{\frac{j-1}{2}} \E_\alpha[v]\big)\big(y \partial_y\bDelta_b^{\frac{h-1}{2}} \phi_\lambda\big)\,.
\end{aligned}
\]
Let us call $\alpha_j :=  \alpha - \frac{j-1}{2}$. Since
$\alpha >[s]=k\geq 2$, $j, h \geq 1$ and $j+h \leq k$, we have
\begin{equation}\label{eq:tech06}
[\alpha_j] \geq k-j+1 \geq 1 \quad\text{and}\quad[\alpha_j -1]\geq k-j \geq 1\,.
\end{equation}
Thus we can use \eqref{eq:Wbdeltam} and $ii)$ in Lemma \ref{L:R1} to infer that there exist some constants $c>0$, not depending on $v$, such that
\begin{equation}\label{eq:tech04}
\begin{aligned}
&\partial_{x_\ell}\bDelta_b^{\frac{j-1}{2}} \E_\alpha[v] = c \E_{\alpha_j}[\tilde v_\ell]\,\\
&y^{-1}\partial_y\bDelta_b^{\frac{j-1}{2}} \E_\alpha[v] = c \E_{\alpha_j-1}[\overline v] \,,
\end{aligned}
\quad\quad \text{ where }\quad\quad 
\begin{aligned}
&\tilde v_\ell = \partial_{x_\ell}(-\Delta)^{\frac{j-1}{2}} v \in C^\infty_c(\R^n)\,,\\
& \overline v = (-\Delta)^{\frac{j+1}{2}}v \in C^\infty_c(\R^n)\,.
\end{aligned}
\end{equation}
On the other hand, by (\ref{eq:lambdaA}) we get
\begin{equation}\label{eq:tech05}
\begin{aligned}
&\partial_{x_\ell}\bDelta_b^{\frac{h-1}{2}}\phi_\lambda = \lambda^{-h} (\tilde \phi_{\ell})_\lambda\\
&y\partial_y\bDelta_b^{\frac{h-1}{2}} \phi_\lambda  = \lambda^{-(h-1)} \overline \phi_\lambda\,,
\end{aligned}
\quad\quad \text{ where }\quad\quad 
\begin{aligned}
&\tilde \phi_\ell = \partial_{x_\ell}\bDelta_b^\frac{h-1}{2}\phi \in \mathcal{A}_0\,,\\
&\overline \phi = y\partial_y\bDelta_b^{\frac{h-1}{2}} \phi \in \mathcal{A}_0\,.
\end{aligned}
\end{equation}
We point out that on the support of $\tilde \phi_\ell$, $\overline \phi \in \mathcal{A}_0$ it holds $\lambda^{-1}\leq c|z|^{-1}$. Therefore
\[ 
\begin{aligned}
|\bnabla_{\!b}^{1+k-(j+h)}\big(\bnabla_{\!b}^j \E_\alpha[v] \cdot \bnabla_{\!b}^h \phi_\lambda \big)|^2 &\leq c\sum_{\ell=1}^n |z|^{-2h}|\bnabla_{\!b}^{1+k-(j+h)}((\tilde \phi_\ell)_\lambda\E_{\alpha_j}[\tilde v_\ell] )|^2 \\
&+c|z|^{-2(h-1)}|\bnabla_{\!b}^{1+k-(j+h)}\big( \overline \phi_\lambda\E_{\alpha_j-1}[\overline v] \big)|^2\,.
\end{aligned}
\]
Next we notice that, by \eqref{eq:tech06} and Lemma \ref{L:extreg1} (see also Remark \ref{R:ck_dense}), we have
\[
\begin{gathered}
(\tilde \phi_\ell)_\lambda\E_{\alpha_j}[\tilde v_\ell]  \in C^{2[\alpha_j]}_{c; \me}(\R^{n+1}) \subset C^{1+k-j}_{c; \me}(\R^{n+1}) \subset \D_\me^{1+k-j;b}(\R^{n+1})\,,\\
 \overline \phi_\lambda\E_{\alpha_j-1}[\overline v]  \in C^{2[\alpha_j-1]}_{c; \me}(\R^{n+1})\subset C^{k-j}_{c; \me}(\R^{n+1}) \subset \D_\me^{k-j;b}(\R^{n+1})\,.
\end{gathered}
\]
As a consequence, we can use Theorem \ref{T:Hardy} to obtain
\[
\begin{aligned}
\iint\limits_{\R^{n+1}}|y|^b|z|^{-2h}|\bnabla_{\!b}^{1+k-(j+h)}\big((\tilde \phi_\ell)_\lambda\E_{\alpha_j}[\tilde v_\ell] \big)|^2 &\leq \sum_{i=0}^{1+k-j}\iint\limits_{\R^{n+1}}|y|^b|z|^{-2i}|\bnabla_{\!b}^{1+k-j-i}\big((\tilde \phi_\ell)_\lambda\E_{\alpha_j}[\tilde v_\ell] \big)|^2 \\
&= \nn (\tilde \phi_\ell)_\lambda\E_{\alpha_j}[\tilde v_\ell] \nn_{1+k-j;b}^2 \leq c\| (\tilde \phi_\ell)_\lambda\E_{\alpha_j}[\tilde v_\ell] \|_{1+k-j;b}^2\,,\\
\iint\limits_{\R^{n+1}}|y|^b|z|^{-2(h-1)}|\bnabla_{\!b}^{1+k-(j+h)}\big( \overline \phi_\lambda\E_{\alpha_j-1}[\overline v]  \big)|^2 &\leq \sum_{i=0}^{k-j}\ \iint\limits_{\R^{n+1}}|y|^b|z|^{-2i}|\bnabla_{\!b}^{k-j-i}\big( \overline \phi_\lambda\E_{\alpha_j-1}[\overline v] \big)|^2 \\
&= \nn  \overline \phi_\lambda \E_{\alpha_j-1}[\overline v] \nn_{k-j;b}^2 \leq c \| \overline \phi_\lambda\E_{\alpha_j-1}[\overline v] \|_{k-j;b}^2\,.
\end{aligned}
\]
Summing up, we proved that
\begin{equation}\label{eq:tech99}
\iint\limits_{\R^{n+1}}|y|^b |\bnabla_{\!b}^{k+1-(j+h)}\big(\bnabla_{\!b}^j \E_\alpha[v] \cdot \bnabla_{\!b}^h \phi_\lambda \big)|^2 \leq c\Big( \|\overline \phi_\lambda \E_{\alpha_j-1}[\overline v] \|_{k-j;b}^2 +\sum_{\ell=1}^n \| (\tilde \phi_\ell)_\lambda\E_{\alpha_j}[\tilde v_\ell] \|_{1+k-j;b}^2\Big)\,,
\end{equation}
for a constant $c>0$ which does not depend on $\lambda$. 

Next, we show that 
\begin{equation}\label{eq:tech08}
\|\overline \phi_\lambda\E_{\alpha_j-1}[\overline v]  \|_{k-j;b}^2 +\sum_{\ell=1}^n \| (\tilde \phi_\ell)_\lambda\E_{\alpha_j}[\tilde v_\ell] \|_{1+k-j;b}^2 \leq c \|\E_\alpha[v] \|_{1+k;b}^2\,.
\end{equation}

Notice that $\overline v, \tilde v_\ell \in C^{\infty}_c(\R^n)$ and $\overline \phi, \tilde \phi_\ell \in \mathcal{A}_0$, see \eqref{eq:tech04} and \eqref{eq:tech05}, respectively. Moreover, \eqref{eq:tech06} holds, thus we are in the position to use  the inductive assumption. 
Taking  in \eqref{eq:crucial2} $v = \tilde v_\ell$, $\phi = \tilde \phi_\ell$ first, and then $v = \overline v$, $\phi = \overline \phi$, we get, respectively, 
\[
\|(\tilde \phi_\ell)_\lambda \E_{\alpha_j}[\tilde v_\ell] \|_{1+k-j;b}^2 \leq c\| \E_{\alpha_j}[\tilde v_\ell]\|_{1+k-j;b}^2~,\quad
\|\overline \phi_\lambda\E_{\alpha_j-1}[\overline v]  \|_{k-j;b}^2 \leq c \|\E_{\alpha_j-1}[\overline v] \|_{k-j;b}^2\,.
\]
Next, \eqref{eq:Wbdeltam} and \eqref{eq:tech04} give $\E_{\alpha_j}[\tilde v_\ell] = c\bnabla_{\!b}^{j-1}\E_\alpha[\partial_{x_\ell}v]$ and $\E_{\alpha_j-1}[\overline v] = c\bnabla_{\!b}^{j+1}\E_\alpha[v]$. Moreover, 
\[
\sum_{\ell=1}^n\| \E_\alpha[\partial_{x_\ell}v]\|^2_{k;b} = c\sum_{\ell=1}^n\irn|(-\Delta)^{\frac{s-1}{2}}\partial_{x_\ell}v|^{2} =c\irn|\Dshalf v|^{2}  = c\|\E_\alpha[v] \|^2_{1+k;b}\,,
\]
by Lemma \ref{L:Walphasumm}. Therefore
\[
\begin{aligned}
\|\overline \phi_\lambda\E_{\alpha_j-1}[\overline v]  \|_{k-j;b}^2 +\sum_{\ell=1}^n \| (\tilde \phi_\ell)_\lambda\E_{\alpha_j}[\tilde v_\ell] \|_{1+k-j;b}^2 &\leq c\big( \|\bnabla_{\!b}^{j+1}\E_\alpha[v] \|^2_{k-j;b}+ \sum_{\ell=1}^n\| \bnabla_{\!b}^{j-1}\E_\alpha[\partial_{x_\ell}v]\|^2_{1+k-j;b}\big)\\
&= c \big( \|\E_\alpha[v] \|^2_{1+k;b} + \sum_{\ell=1}^n\| \E_\alpha[\partial_{x_\ell}v]\|^2_{k;b}\big) = c\|\E_\alpha[v] \|_{1+k;b}^2\,,
\end{aligned}
\]
which proves \eqref{eq:tech08}. By \eqref{eq:tech99} and \eqref{eq:tech08} we obtain \eqref{eq:tech03}, and the proof of  (\ref{eq:crucial2}) is complete.

By taking $v = u$ and $\phi = \f$ in \eqref{eq:crucial2} we  get \eqref{eq:crucial}. This concludes the proof of $i)$ in Lemma \ref{L:semplice}.
\QED

\paragraph{Acknowledgments.}
This work is partially supported by the PRID  Projects {\sc priden} and {\sc vaproge}, 
 Universit\`a di Udine.

\end{document}